\documentclass[11pt,thmsa,emstex]{article}
\usepackage[left=3cm,right=3cm,top=3cm,bottom=3cm]{geometry}
\usepackage{amsfonts}
\usepackage{t1enc}
\usepackage{amssymb}
\usepackage{epsfig}
\usepackage[french,english]{babel}
\usepackage{amsmath}
\usepackage{graphics}
\usepackage{layout}
\usepackage{latexsym}
\usepackage{color}
\DeclareMathOperator{\Tr}{Tr}

\def\mylabel#1{\label{#1}}

\newtheorem{theorem}{Theorem}[section]
\newtheorem{lemma}[theorem]{Lemma}
\newtheorem{corollary}[theorem]{Corollary}
\newtheorem{proposition}[theorem]{Proposition}
\newtheorem{example}[theorem]{Example}
\newtheorem{remark}[theorem]{Remark}

\newtheorem{hypothesis}[theorem]{Hypothesis}

\newtheorem{exercice}[theorem]{Exercice}

\def\bit{\begin{itemize}}

\def\eit{\end{itemize}}

\reversemarginpar   

\def\bc{\begin{center}}

\def\ec{\end{center}}

\def\bthm{\begin{theorem}}

\def\ethm{\end{theorem}}

\def\bcor{\begin{corollary}}

\def\ecor{\end{corollary}}

\def\bprop{\begin{proposition}}

\def\eprop{\end{proposition}}

\def\blem{\begin{lemma}}

\def\elem{\end{lemma}}

\def\bex{\begin{example}}

\def\eex{\end{example}}

\def\bexo{\begin{exercice} \rm }

\def\eexo{\end{exercice} }

\def\brem{\begin{remark}}

\def\erem{\end{remark}$\Box$}

\def\prf{{\bf Proof }}

\def\bdes{\begin{description}}

\def\edes{\end{description}}

\def\ita{\item[(a)]}

\def\itb{\item[(b)]}

\def\itc{\item[(c)]}

\def\itd{\item[(d)]}

\def\iti{\item[(i)]}

\def\itii{\item[(ii)]}

\def\itiii{\item[(iii)]}

\def\itiv{\item[(iv)]}

\def\beq{\begin{equation}}

\def\eeq{\end{equation}}

\def\ben{\begin{enumerate}}

\def\een{\end{enumerate}}

\def\beqar{\begin{eqnarray}}

\def\eeqar{\end{eqnarray}}

\def\beqarr{\begin{eqnarray*}}

\def\eeqarr{\end{eqnarray*}}

\def\qed{\hspace{.1in}{\bf QED}\\[2ex]}

\def\prf{{\bf Proof }\hspace{.1in}}

\def\L{{\cal L}}
\def\D{{\cal D}}

\def\Pr{{\mathsf P}}
\def\E{{\mathsf E}}

\def\B{{\cal B}}
\def\P{{\mathbb P}}

\def\Ind{{\mathbf 1}}
\def\RR{{\mathbb R}}  
\def\Rp{{\mathbb R}_+}   

\def\NN{{\mathbb N}}

\def\XX{{\mathcal X}}
\def\rar{\rightarrow}

\def\eps{\varepsilon}

\def\1{{\rm 1\mskip-4.4mu l}}
\begin{document}
\title{Random Switching between Vector Fields Having a Common Zero}
\author{Michel Bena\"{\i}m, Edouard Strickler\\ Institut de Math\'ematiques\\Universit\'e de Neuch\^atel, Switzerland}
\maketitle
\begin{abstract}
Let $E$ be a finite set, $\{F^i\}_{i \in E}$ a family of vector fields on $\RR^d$ leaving positively invariant a compact set $M$ and having a common zero $p \in M.$ We consider a piecewise deterministic Markov process $(X,I)$ on $M \times E$ defined by $\dot{X}_t = F^{I_t}(X_t)$ where $I$ is a jump process controlled by $X:$ $\Pr(I_{t+s} = j | (X_u, I_u)_{u \leq t}) = a_{i j}(X_t) s + o(s)$ for $i \neq j$ on $\{I_t = i \}.$

We show that the behavior of $(X,I)$ is mainly determined by the behavior of the linearized process $(Y,J)$ where $\dot{Y}_t = A^{J_t} Y_t,$
$A^i$ is the Jacobian matrix of $F^i$ at $p$ and $J$ is the jump process with rates $(a_{ij}(p)).$
 We introduce two quantities $\Lambda^-$ and $\Lambda^+$  respectively
defined as the {\em minimal} (respectively {\em maximal}) {\em growth rate} of $\|Y_t\|,$ where the minimum (respectively maximum) is taken over all the ergodic  measures of the angular process $(\Theta, J)$ with $\Theta_t = \frac{Y_t}{\|Y_t\|}.$ It is shown that $\Lambda^+$ coincides with the top Lyapunov exponent (in the sense of ergodic theory) of $(Y,J)$ and that under general assumptions $\Lambda^- = \Lambda^+.$ We then prove that, under certain irreducibility conditions, $X_t \rar p$ exponentially fast when $\Lambda^+ < 0$ and $(X,I)$ converges in distribution at an exponential rate toward a (unique) invariant measure supported by $M \setminus  \{p\} \times E$ when $\Lambda^- > 0.$
Some applications to certain epidemic models in a fluctuating environment are discussed and illustrate our results.
\end{abstract}
\newpage
\tableofcontents
\newpage

\section{Introduction}
\mylabel{sec:intro}

Let $E$ be a finite set and  $\mathrm{F}  = \{F^i\}_{i \in E}$ a family of $C^2$ globally integrable vector fields on $\RR^d.$
For each $i \in E$ we let $\Psi^i = \{\Psi^i_t\}$   denote the  flow induced by $F^i.$
We assume  throughout that there exists a  closed set $M \subset \RR^d$  which is {\em positively invariant} under each $\Psi^i.$ That is $$\Psi^i_t(M) \subset M$$ for all $t \geq 0.$

Consider a Markov process $Z = (Z_t)_{t \geq 0}, Z_t = (X_t, I_t),$  living on $M \times E$ whose infinitesimal generator acts on functions $g : M \times E \mapsto \RR,$ smooth in the first variable, according to the formula
\beq
\label{eq:defcL}
{\cal L}g(x,i) =  \langle F^i(x), \nabla g^i(x) \rangle + \sum_{j \in E} a_{ij}(x) (g^j(x) - g^i(x)),
\eeq
where $g^i(x)$ stands for $g(x,i)$ and  $a(x) = (a_{ij}(x))_{i,j \in E}$ is an irreducible {\em rate} matrix continuous in $x.$ Here, by a rate matrix, we mean a matrix  having nonnegative off diagonal entries and zero diagonal entries.

In other words, the dynamics of $X$ is given by an ordinary differential equation
\beq
 \frac{dX_t}{dt} = F^{I_t}(X_t),
\label{eq:pdmp}
\eeq
while  $I$ is a continuous time jump process taking values in $E$ controlled by $X:$
$$\Pr(I_{t+s} = j | {\cal F}_t, I_t = i)  = a_{ij}(X_t) s + o(s) \mbox{ for } j \neq i \mbox{ on } \{I_t = i\},$$
where ${\cal F}_t = \sigma ((X_s,I_s) \: : s \leq t\}.$

This class of processes belongs to the wider class of {\em Piecewise Deterministic Markov Processes} (PDMPs), a term coined by Davis \cite{Dav84},
and has recently been the focus of much attention. Criteria, based on irreducibility and H{\"o}rmander type conditions,
 ensuring  uniqueness and absolute continuity of an invariant probability measure  have been obtained by Bakhtin and Hurth  \cite{bakhtin&hurt} for constant jump rates ($a_{ij}(x) = a_{ij})$ and by Bena{\"i}m, Le Borgne, Malrieu and Zitt   \cite{BMZIHP} for more general rates. Exponential convergence (in total variation) toward this measure and a support theorem, describing the support of the law of $(Z_t)_{z \geq 0}$  are also proved in \cite{BMZIHP} when $M$ is compact (see also \cite{BHS18}). In the one dimensional case (i.e $d = 1$) smoothness properties of the invariant measure are thoroughly investigated by Bakhtin, Hurth and Mattingly   \cite{bakhtin&hurt&matt}. When irreducibility fails to hold, the support of invariant probabilities can be determined in terms of invariant control sets of an associated deterministic control system (see Bena{\"i}m, Colonius and Lettau \cite{BCL16}).  When the vector fields are exponentially asymptotically stable in "average", exponential convergence toward an invariant measure are obtained for  Wassertein distances   by Bena{\"i}m, Le Borgne, Malrieu and Zitt \cite{ecp},    Cloez and Hairer \cite{CH13}. Several examples, either linear (Bena{\"i}m, Le Borgne, Malrieu and Zitt \cite{BLMZexample}, Lawley, Mattingly and Reed \cite{Lawley&matt&Reed}, Lagasquie \cite{lagasquie16}),  or nonlinear (Bena{\"i}m and Lobry \cite{BL16}, Malrieu and Hoa Phu \cite{MalHoa16}) show that the behavior of the process is not solely determined by the dynamics of the $\Psi^i$ but can be highly sensitive to the switching rates. We refer the reader to the  recent overview by Malrieu \cite{Mal15},  describing these results among others.
\medskip

In the present paper we will investigate the behavior of the process $Z$ under the following two conditions:

\bdes
\item[C1] The origin lies in $M$ and is a common equilibrium:  $$F^i(0) = 0 \mbox{ for all } i \in E;$$
\item[C2] The set $M$ is compact and {\em locally star shaped at the origin}, meaning that there exists
$\delta > 0$ such that
$$x \in M \mbox{ and } \|x\| \leq \delta \Rightarrow [0,x] \subset M,$$ where $[0,x] = \{ tx, \: t \in [0,1]\}$.
\edes
 Compactness of $M$ is assumed here  for simplicity, but some of the (local) results generalise to noncompact sets. The global results can be extended  provided we can control the behaviour of the process near infinity, for instance with a suitable Lyapunov function (see Section \ref{sec:noncompact}).

 Briefly put, our main result is that the long term behavior of the process is determined by the behavior of the process obtained by linearization at the origin and,
  under suitable irreducibility and hypoellipticity conditions, by the top Lyapunov exponent of the linearized system.
  If negative, then $X = (X_t)$ converges almost surely and exponentially fast to zero. If positive, and $X_0 \neq 0,$
  the empirical occupation measure (respectively  the law) of $Z$ converge almost surely  (respectively in total variation at an exponential rate)
  toward a unique  probability measure putting zero mass on $\{0\} \times E.$ Such a  correspondence between the sign  of the top Lyapunov exponent and the behavior of nonlinear system is reminiscent of the results obtained by Baxendale  \cite{Bax90} and others    for Stratonovich stochastic differential equations (see  \cite{Bax90} and the references therein, and Hening, Nguyen and Yin \cite{Hen16} for similar recent results in the context of population dynamics).

  Our proofs rely,  on one hand, on the qualitative theory of PDMPs (as developed in \cite{bakhtin&hurt} and \cite{BMZIHP}) and, on the other hand, on some recent results on {\em stochastic persistence} (Bena{\"i}m \cite{Ben14}) strongly inspired by the seminal works of Schreiber, Hofbauer and their co-authors on {\em persistence}, first developed for purely deterministic systems (Schreiber \cite{S00}, Garay and Hofbauer \cite{GH03}, Hofbauer and Schreiber \cite{HofSch04}) and later for certain stochastic systems (Bena\"{\i}m, Hofbauer  and Sandholm \cite{BHW}, Bena{\"i}m and  Schreiber \cite{BS09},  Schreiber, Bena\"{\i}m and Atchade \cite{SBA11}, Schreiber \cite{Sch12}, Roth and Schreiber \cite{SR14}).

Our original motivation was to analyze the behavior of certain epidemic models evolving in a fluctuating environment.
 A famous, and now classical, deterministic model of infection is given by the Lajmanovich and Yorke differential equation (\cite{LajYorke}). This equation  leaves positively invariant the unit cube of $\RR^d$ and models the evolution of the infection level between $d$ groups. Depending on the parameters of the model (the environment), either the disease dies out (i.e~ all the  trajectories converge to the origin) or stabilizes (i.e~ all non zero trajectories converge toward a unique positive equilibrium). Deterministic switching between several environment have been recently considered by Ait Rami, Bokharaie, Mason and Wirth  \cite{Masonetal14}. The results here allow to describe the behavior of the process when switching between environment evolves randomly. In particular we can  produce paradoxical examples for which, although each deterministic dynamics leads to the extinction (respectively persistence) of the disease, the random switching process leads to persistence (respectively extinction) of the disease.

 \subsection{Outline of contents}
Section \ref{sec:linear} considers the linearized system $(Y,J)$ where $\dot{Y}_t = A^{J_t} Y_t, A^i = DF^i(0)$ (the Jacobian of $F^i$ at $0$) and $J$ is the jump process with rate matrix $(a_{ij}) = (a_{ij}(0)).$
 We introduce two quantities $\Lambda^-$ and $\Lambda^+$  respectively
defined as the {\em minimal} (respectively {\em maximal}) {\em growth rate} of $\|Y_t\|,$ where the minimum (respectively maximum) is taken over all the ergodic  measures of the angular Markov process $(\Theta, J)$ with $\Theta_t = \frac{Y_t}{\|Y_t\|}.$ It is shown (Proposition \ref{th:lyapou}) that $\Lambda^+$ coincides with the top  Lyapunov exponent (in the sense of ergodic theory) of $(Y,J)$ and some conditions are given ensuring that $\Lambda^- = \Lambda^+,$ first for arbitrary $A^i$s (Proposition \ref{th:uniqueLambda}) and then for Metzler matrices  (Proposition \ref{th:Perron}).

The main results of the paper are stated in Section \ref{sec:main}.
 \begin{itemize}
 \item If $\Lambda^+ < 0,$ $X_t \rar 0$ exponentially fast, locally (i.e~ for $\|X_0\|$ small enough), with positive probability. If furthermore $0$ is accessible, convergence is global and almost sure (Theorem \ref{th:extinct}).
       \item If $\Lambda^- > 0$ and $X_0 \neq 0$, the process is {\em persistent} in the sense that  weak limit points of its empirical occupation measure are almost surely invariant probabilities over $M \setminus \{0\} \times E$ (Theorem \ref{th:persist1}). If in addition the $F^i$s satisfy a certain H{\"o}rmander-type bracket condition at some accessible point, then there is a unique invariant probability on $M \setminus \{0\} \times E$ toward which the empirical occupation measure converges almost surely (Theorem \ref{th:persist2}). Under a strengthening of the bracket condition, the distribution of the process converges also exponentially fast in total variation (Theorem \ref{th:persist3}).
 \end{itemize}

Section \ref{sec:epidemic} discusses some applications of our results to certain epidemic models in a fluctuating environment. The focus is on the situation where the $F^i$s are given by Lajmanovich and Yorke type vector fields \cite{LajYorke} (or more generally sub homogeneous cooperative systems in the sense of Hirsch \cite{Hirsch94}). Several examples are analyzed and a theorem proving exponential convergence of the distribution (for a certain Wasserstein distance) in absence of the bracket condition is stated (Theorem \ref{th:cvexpo}).

Sections \ref{sec:mainproofs} and \ref{sec:proofcvexpo}  are devoted to the proofs of Theorems \ref{th:extinct}, \ref{th:persist1}, \ref{th:persist2}, \ref{th:persist3} and \ref{th:cvexpo}. The proofs of certain results stated in Section \ref{sec:linear} are given in  appendix (Section \ref{sec:appendix}) for convenience.

 \subsection{Notation}
The following notation will be used throughout: $\langle \cdot, \cdot  \rangle$ denotes the Euclidean inner product in $\RR^d,$  $\| \cdot \|$ the associated norm,
 $B(x, r) = \{y \in \RR^d \: : \|y -x\| \leq r\}$
 the closed ball centered at $x$ with radius $r$ and $S^{d-1} = \{x \in \RR^d \: : \|x\| = 1\}$ the unit sphere.

\paragraph{Notation for Markov processes} For any polish space $\mathcal{X}$ such as $M, S^{d-1}, E,$ $M \times E,$  equipped with its Borel sigma-field, we let ${\cal P}(\mathcal{X})$
  denote the set of (Borel) probabilities over $\mathcal{X}.$ We shall consider below certain Markov processes  $\tilde{Z}$ (like $Z$)  taking values in $\mathcal{X}$ with {\em cad-lag} (right continuous, left limit) paths. Given such a process and $\mu \in {\cal P}(\mathcal{X})$ we let $\mathbb{P}^{\tilde Z}_{\mu}$ denote the law of $\tilde{Z}$
   on the Skorokhod space $D(\RR^+, \mathcal{X})$ when $\tilde{Z}_0$ has law $\mu.$  As usual,   $\mathbb{P}^{\tilde Z}_{z}$ stands for $\mathbb{P}^{\tilde Z}_{\delta_{z}}$ for all $z \in \mathcal{X}.$ The {\em Markov semi-group} induced by $\tilde{Z},$ denoted $(P_t^{\tilde{Z}})_{t \geq 0},$ acts on bounded measurable functions $f : \mathcal{X} \mapsto \RR$ according to the formula $$P_t^{\tilde{Z}} f(z) = \mathbb{E}_z(f(\tilde{Z}_t)) = \int f(\eta(t)) d \mathbb{P}^{\tilde Z}_{z}(\eta).$$ By duality it acts on ${\cal P}(\mathcal{X})$ by
   $$(\mu P_t^{\tilde{Z}}) f = \mu (P_t^{\tilde{Z}} f),$$ where here and throughout $\mu f$ stands for $\int f d\mu.$
Probability $\mu \in {\cal P}(\mathcal{X})$ is said {\em invariant} for $\tilde{Z}$ provided $\mu P_t^{\tilde{Z}} = \mu$ for all $t \geq 0.$ It is called {\em ergodic} if, in addition of being invariant, the only  bounded measurable functions $f : \mathcal{X} \mapsto \RR $ for which $\sup_{t \geq 0} \mu( |P_t^{\tilde{Z}} f - f|) = 0$ are $\mu$-almost surely constant.

   We let ${\cal P}_{inv}^{\tilde{Z}} \subset {\cal P}(\mathcal{X})$ denote the (possibly empty) set of {\em invariant probabilities} of $\tilde{Z}$ and ${\cal P}_{erg}^{\tilde{Z}} \subset {\cal P}_{inv}^{\tilde{Z}}$  the  subset of {\em ergodic probabilities}. Recall that ${\cal P}_{erg}^{\tilde{Z}}$ can also be defined as the  set of extremal points of  ${\cal P}_{inv}^{\tilde{Z}}.$

   A key property, that will be used later without further notice, is that whenever $\mu \in {\cal P}_{inv}^{\tilde{Z}}$  (respectively
    $\mu \in {\cal P}_{erg}^{\tilde{Z}}$),  $\mathbb{P}^{\tilde Z}_{\mu}$ is invariant (respectively ergodic), in the sense of ergodic theory, for  the shift $\mathbf{\Theta} = (\mathbf{\Theta}_t)_{t \geq 0}$ on  $D(\RR^+, \mathcal{X});$ where $$\mathbf{\Theta}_t(\eta)(s) = \eta(t+s).$$
    We refer the reader to Meyn and Tweedie (\cite{MT1}, chapter 17) for a proof and more details.
    \paragraph{Accessibility} Let $\tilde{\mathrm{F}} = \{\tilde{F}^i\}_{i \in E}$ be a family of bounded $C^1$ vector fields on $\RR^d$ indexed by $E.$ For instance $\tilde{\mathrm{F}} = \mathrm{F}.$
    We let $\mathsf{co}(\tilde{\mathrm{F}})$ denote the compact convex set valued mapping defined by
    $$\mathsf{co}(\tilde{\mathrm{F}})(x) = \{ \sum_{j \in E} \alpha_j \tilde{F}^j(x) \: : \alpha_j \geq 0, \sum_{j \in E} \alpha_j = 1\}.$$
    Given a closed set $A \subset  \RR^d$ and $B \subset \RR^d$ we say that $A$ is {\em $\tilde{F}$-accessible} from $B$ if for every neighborhood $U$ of $A$
    and every $x \in B,$ there exists a (absolutely continuous) function $\eta :  \Rp \mapsto \RR^d,$ solution to the differential inclusion
    $$\left \{\begin{array}{l}
\dot{\eta}(t) \subset \mathsf{co}(\tilde{\mathrm{F}})(\eta(t)) \\
\eta(0) = x\end{array} \right.$$  such that $\eta(t) \in U$ for some $t > 0$.
 An equivalent formulation (see e.g~ Theorem 2.2 in \cite{BCL16}) is that $A$ is {\em reachable} from $B$ by the control system
 $$\left \{ \begin{array}{l} \dot{x} = \sum_j \tilde{F^j} (x) v_j(t) \\ x(0) = x \end{array} \right.$$ where the control $v \in D(\Rp, \{e_j\}_{j \in E})$ with $\{e_j\}_{j \in E}$  the canonical basis of $\RR^E.$  Note that this notion is what is called $D$-approachability in \cite{bakhtin&hurt}.

\section{The Linearized system}
\label{sec:linear}
Let, for $i \in E, A^i = DF^i(0)$ denote the Jacobian matrix of $F^i$ at the origin.
We let $C_M \subset \RR^d$ denote the cone defined as
$$C_M = \overline{\{t x : \: t \geq 0, x \in M, \|x\| \leq \delta\}}$$
where $\delta$ is like in condition $C2.$ Here, $\overline{B}$ stands for the closure of $B$.
\brem
One can check that the definition of $C_M$ does not depend on the choice of $\delta$, provided $\delta$ satisfies condition $C2$.
\erem
\blem
For all $t \geq 0$ $e^{t A^i} C_M \subset C_M.$
\elem
\prf We set $D_M = \{t x : \: t \geq 0, x \in M, \|x\| \leq \delta\}$ and first prove that $e^{t A^i} D_M \subset C_M.$ The lemma will be then induced by continuity of $e^{t A^i}$.  Let $x \in D_M.$ For $\eps$ small enough, by definition of $D_M$ and continuity of $\Psi^i_t$ at $0$
$\Psi^i_t(\eps x) \in M \cap B(0, \delta).$ Hence $\frac{\Psi^i_t(\eps x)}{\eps} \in C_M$ and letting $\eps \rar 0$ this shows that
$D\Psi^i_t(0) x = e^{t A^i} x \in C_M.$
\qed

Define the {\em linearized system of $Z$ at the origin} as the "linear" PDMP $(Y,J)$ living on $C_M \times E$ whose generator $L$ is given by
$$L g(y,i) = \langle A^i y, \nabla g^i(y) \rangle + \sum_{j \in E} a_{ij} (g^j(y) - g^i(y)),$$
where $$a_{ij} = a_{ij}(0).$$
A trajectory $(Y_t,J_t)_{t \geq 0}$ with initial condition $(y,i)$ is then obtained as a solution to
\beq
\label{dY}
 \left \{
\begin{array}{l}
  \frac{dY_t}{dt} = A^{J_t} Y_t \\
  Y_0 = y, \\
\end{array} \right .\eeq
where
$(J_t)$ is a continuous time Markov process on $E$ with jump rates
$(a_{ij})$ based at $J_0 = i.$

By irreducibility of $(a_{ij}), J$ has a unique invariant probability $p = (p_i)_{i \in E},$ characterized by
$$\forall i \in E, \, \sum_j (p_j a_{ji} - p_i  a_{ij}) = 0.$$
Whenever $y \neq 0$ the  {\em polar decomposition} $$(\Theta_t = \frac{Y_t}{\|Y_t\|}, \rho_t = \|Y_t\|) \in S^{d-1} \cap C_M \times \Rp$$ is well defined and
(\ref{dY})  can be rewritten as
\beq
\label{dThetarho}
 \left \{
\begin{array}{l}
  \frac{d\Theta_t}{dt} = G^{J_t}(\Theta_t) \\
  \frac{d \rho_t}{dt} =  \langle A^{J_t} \Theta_t, \Theta_t \rangle \rho_t,
\end{array}
\right.
\eeq
where for
 all $i \in E$ $G^i$ is the vector field on $S^{d-1}$ defined by
\beq
\label{eq:defGi}
G^i(\theta) = A^i \theta - \langle A^i \theta, \theta \rangle \theta.
\eeq
\brem 
For stochastic differential equations, the idea of introducing, this polar decomposition goes back to Hasminskii \cite{Has60} and has proved to be a fundamental tool for analyzing linear stochastic differential equations (see e.g~ \cite{Bax90}), linear random dynamical systems (see e.g~chapter 6 of Arnold \cite{Arn98}) and more recently certain linear PDMPs in \cite{BLMZexample},~\cite{Lawley&matt&Reed} or~\cite{lagasquie16}.
\erem

With obvious notation, the processes $$(\Theta, \rho, J) = ((\Theta_t, \rho_t, J_t))$$ and $$(\Theta, J) = ((\Theta_t, J_t))$$ are two PDMPs respectively living on $S^{d-1} \cap C_M \times \Rp \times E$ and $S^{d-1} \cap C_M \times E.$

By compactness of $S^{d-1} \cap C_M$ and Feller continuity of $(\Theta, J)$ (see \cite{BMZIHP}, Proposition 2.1),  ${\cal P}_{inv}^{(\Theta, J)}$ is a nonempty compact (for the topology of weak* convergence) subset of ${\cal P}(S^{d-1} \cap C_M \times E).$

\subsection{Average growth rates}
Define, for each $\mu \in {\cal P}_{inv}^{(\Theta, J)},$ the $\mu$-{\em average growth rate} as
\beq
\label{deflambdamu}
\Lambda(\mu) = \int \langle A^i \theta, \theta \rangle \mu(d\theta di) = \sum_{i \in E} \int_{S^{d-1} \cap C_M} \langle A^i \theta, \theta \rangle \mu^i(d\theta),
\eeq where
$\mu^i(.)$ is the measure on $S^{d-1} \cap C_M$ defined by $$\mu^i(A) = \mu (A \times \{i\}).$$
Note that when $\mu$ is ergodic, by equation (\ref{dThetarho}) and Birkhoff ergodic theorem
$$\lim_{t \rar \infty} \frac{\log(\rho_t)}{t} = \Lambda(\mu)$$
$\P_{\mu}^{(\Theta, J)}$ almost surely.

Define similarly the {\em extremal average growth rates} as the numbers
\beq
\label{deflambdamumax}
\Lambda^{-} = \inf\{ \Lambda(\mu) \: : \mu \in {\cal P}_{erg}^{(\Theta, J)} \} \mbox{ and }
\Lambda^+ = \sup \{\Lambda(\mu) \: : \mu \in {\cal P}_{erg}^{(\Theta, J)}\}.
\eeq
The following rough estimate is a direct consequence of (\ref{deflambdamu}). Recall that $p = (p_i)_{i \in E}$ is the invariant probability of $J.$
\blem
\label{lem:roughestimate}
$$\sum_i p_i \lambda_{min} (\frac{A^i + (A^i)^T}{2}) \leq \Lambda^{-} \leq \Lambda^{+} \leq
\sum_i p_i \lambda_{max} (\frac{A^i + (A^i)^T}{2}),$$
where $\lambda_{min}$ (respectively $\lambda_{max}$) denotes the smallest (respectively largest) eigenvalue.
\elem
The signs of  $\Lambda^-$ and $\Lambda^+$ will play a crucial role for determining the asymptotic behavior of the non linear process $Z.$
But before stating our main results, it is interesting to compare them with the usual Lyapunov exponents given by the multiplicative ergodic theorem.
\subsection{Relation with Lyapunov exponents}
\label{sec:lyapou}
Set $\Omega = D(\Rp, E)$ and
for $\omega \in \Omega$ and $y \in \RR^d,$ let $$t \mapsto \varphi(t,\omega) y$$ denote  the solution to the linear  differential equation
$$\dot{y} = A^{\omega_t} y$$ with initial condition $\varphi(0,\omega) y = y.$

Then,
$\varphi$ is a {\em linear random dynamical system} over the ergodic dynamical system
$(\Omega,  \mathbb{P}^J_{p}, \mathbf{\Theta}),$ for which the assumptions of the multiplicative ergodic theorem are easily seen to be satisfied (see e.g~ \cite{Arn98}, Theorem 3.4.1 or Colonius and Mazanti \cite{ColMa15}). Thus, according to this theorem, there exist $1 \leq \tilde{d} \leq d,$ numbers
$$\lambda_{\tilde{d}} < \ldots < \lambda_1,$$ called {\em the Lyapunov exponents} of $(\varphi, \mathbf{\Theta})$ , a Borel set $\tilde{\Omega}  \subset \Omega$ with $\mathbb{P}^J_{p}(\tilde{\Omega}) = 1,$ and for each $\omega \in \tilde{\Omega}$ distinct vector spaces
$$\{0\} = V_{\tilde{d}+1}(\omega) \subset V_{\tilde{d}}(\omega) \subset \ldots \subset V_{i}(\omega) \ldots \subset V_1(\omega) = \RR^d$$ (measurable in $\omega$) such that
\beq
\label{eq:deflambdai}
\lim_{t \rar \infty} \frac{1}{t} \log \|\varphi(t,\omega) y\| = \lambda_i
\eeq for all $y \in V_i(\omega) \setminus V_{i+1}(\omega)$.
\bprop
\label{th:lyapou}
For all $\mu \in {\cal P}_{erg}^{(\Theta, J)}$
$$\Lambda(\mu) \in \{\lambda_{\tilde{d}}, \ldots, \lambda_1\}.
$$
If furthermore $C_M$ has non empty interior, then $$\Lambda^+ = \lambda_1.$$
\eprop

\brem
{\em The second part of the proposition has already been proven by Crauel \cite[Theorem 2.1 and Corollary 2.2]{C84} in a more general setting. We adapt the arguments of his proof for our specific case.}\erem

\prf
Let $\mu \in   {\cal P}_{erg}^{(\Theta, J)}.$  Then, $\P_{\mu}^{(\Theta,J)}$ almost surely
$$\lim_{t \rar \infty} \frac{1}{t} \log (\|\varphi(t,J) \Theta_0)\|) =
\lim_{t \rar \infty} \frac{1}{t} \int_0^t \langle A^{J_s} \Theta_s, \Theta_s \rangle ds = \Lambda(\mu)$$
The first equality follows from  (\ref{dY}), (\ref{dThetarho}) and the definition of $\varphi(t,\omega).$ The second follows from Birkhoff ergodic theorem.
Therefore, there exists a Borel set $\B \subset  (S^{d-1} \cap C_M) \times \Omega$
such that for all $(\theta,\omega) \in \B$
\beq
\label{eq:lambdaLambda}
\lim_{t \rar \infty} \frac{1}{t} \log (\|\varphi(t,\omega) \theta\|)  = \Lambda(\mu)
\eeq  and $\P_{\mu}^{(\Theta_0,J)}(\B) = 1,$
where $\P_{\mu}^{(\Theta_0,J)}(d\theta d \omega) = \sum_{i \in E} \P_{i}^{J}(d \omega) \mu^i(d\theta)$ is the law of $(\Theta_0, J)$  under $\P_{\mu}^{(\Theta,J)}.$

 Let $\tilde{\Omega} \subset \Omega$ be the set given by the multiplicative ergodic theorem and
$\tilde{\B} = \{(\theta,\omega) \in \B \: : \omega \in \tilde{\Omega}\}.$ Then
 $\P_{\mu}^{(\Theta_0,J)}(S^{d-1} \cap C_M \times \tilde{\Omega}) = \P_{\mu}^{J}(\tilde{\Omega}) = 1.$
Hence $\P_{\mu}^{(\Theta_0,J)}(\tilde{\B}) = 1$  and for all $(\theta, \omega) \in \tilde{\B}$ the left hand side of equality (\ref{eq:lambdaLambda}) equals $\lambda_i$ for some $i.$

It remains to show that $\lambda_1 = \Lambda^+$. For every $\omega$ in the set $\tilde{\Omega}$ given by the multiplicative ergodic theorem, and for all $\theta \in S^{d-1} \cap C_M$, define $$\lambda(\theta,\omega) = \lim_{ t \to \infty}\frac{1}{t} \log (\|\varphi(t,\omega) \theta\|)=\lim_{ t \to \infty}\frac{1}{t}\int_0^t \langle A^{\omega_s} \Theta_s^{\theta}(\omega),\Theta_s^{\theta}(\omega)\rangle ds,$$
where $$\Theta_t^{\theta}(\omega) = \frac{\varphi(t,\omega) \theta}{\|\varphi(t,\omega) \theta \|}.$$  By \eqref{eq:deflambdai}, we have $\lambda(\theta,\omega) = \lambda_1$ for all $\theta \in V_1(\omega) \setminus V_2(\omega) \cap S^{d-1} \cap C_M$. Let $\nu$ denote the normalised Lebesgue measure on $ S^{d-1} \cap C_M$. Because $V_2(\omega)$ is at most an hyperplane and $C_M$ has non empty interior, we get that $\int \lambda(\theta, \omega) d\nu(\theta) = \lambda_1$ for all $\omega \in \tilde{\Omega}$. In particular,
\beq 
\label{eq:intintlambda}
\int_{\Omega} \int_{S^{d-1} \cap C_M} \lambda(\theta, \omega) d\nu(\theta) d\mathbb{P}^J_p(\omega) = \lambda_1.
\eeq Moreover, because $|\langle A^i \theta, \theta \rangle| \leq \max \|A^i\|$, dominated convergence and \eqref{eq:intintlambda} imply that 
\beq
\label{eq:lambdalimint}
\lambda_1 = \lim_{ t \to \infty} \frac{1}{t} \int_{\Omega} \int_{S^{d-1} \cap C_M} \int_0^t \langle A^{\omega_s} \Theta_s^{\theta}(\omega),\Theta_s^{\theta}(\omega) \rangle ds  d\nu(\theta) d\mathbb{P}^J_p(\omega)
\eeq
Now for all $t > 0$, define the probability on $ S^{d-1} \cap C_M \times E$ 
\beq
\label{eq:mut}
\mu_t = \frac{1}{t} \int_0^t (\nu \otimes p)P_s^{(\Theta,J)}ds.
\eeq
By compactness of $ S^{d-1} \cap C_M \times E$, $(\mu_t)_{t\geq 0}$ is tight, and by Feller property of $(\Theta, J)$, every weak limit points of $\mu_t$ belongs to $\mathcal{P}_{inv}^{(\Theta,J)}(S^{d-1} \cap C_M \times E)$. Let $\mu$ be such a limit point, and $(t_n)$ such that $\mu_{t_n} \rightarrow \mu$. Setting $f(\theta,i) = \langle A^i \theta, \theta \rangle$, one has $\mu_{t_n}f \to \mu f = \Lambda(\mu)$. Now \eqref{eq:lambdaLambda}, \eqref{eq:lambdalimint} and Fubini Theorem imply that $\lambda_1 = \lim \mu_{t_n} f = \Lambda(\mu)$, which concludes the proof.
\qed

In the multiplicative ergodic theorem, each Lyapunov exponent $\lambda_i$ comes with an integer $d_i \geq 1$ called its {\em multiplicity} and such that $\sum_{i = 1}^{\tilde{d}} d_i = d$ (see Chapter 3 of \cite{Arn98} for more details).
A consequence of Proposition \ref{th:lyapou} is the following inequality which provides, in some cases,  a simple way to prove that  $\Lambda^+ > 0$, which is often  a sufficient condition to ensure positive recurrence of $Z$ on $M \setminus \{0\} \times E$ (see Propostions \ref{th:uniqueLambda} and \ref{th:Perron} and Theorems \ref{th:persist1} and \ref{th:persist2}).
\bcor
\label{cor:trace}
$$\sum_{i \in E} p_i \Tr(A^i) = \sum_{i = 1}^{\tilde{d}} d_i \lambda_i \leq  d \Lambda^+.$$
\ecor
\prf
 By Jacobi's formula
$$\frac{\log(\det(\varphi(t,\omega)))}{t} = \frac{\int_0^t \Tr(A^{\omega_s})ds}{t} .$$  By Birkhoff ergodic Theorem, the right hand side of this equality converges,  $\mathbb{P}^J_{p}$ almost surely, as $t \rar \infty,$ toward $\sum_i p_i \Tr(A^i);$ and a  by product of the multiplicative ergodic theorem (see e.g~ \cite{Arn98}, Chapter 3, Corollary 3.3.4) is that the left-hand side converges $\mathbb{P}^J_{p}$ almost surely,
as $t \rar \infty,$ toward $\sum_{i = 1}^{\tilde{d}} d_i \lambda_i.$
\qed
\brem
\label{rem:Mierc}
 {\em If the  matrices $A^i$ are {\em Metzler}, meaning that they have off diagonal nonnegative entries, a result due to   Mierczy{\'n}ski (\cite{Mier13}, Theorem 1.3)
 allows to improve the lower bound given in Corollary \ref{cor:trace} We will use this estimate in section \ref{sec:epidemic}, example \ref{ex:nobra}.}
 \erem
\brem
{\em Note that in general  $$\Lambda^- \neq \lambda_{\tilde{d}}.$$
 Here is a simple example based on \cite{BLMZexample}. Assume $E = \{1,2\}$ and $d = 2$ (so that the matrices here are $2 \times 2$). Let $A^1, A^2$ be  $2$ real matrices having eigenvalues with negative real parts and such  that for some $0 < t < 1,$ the eigenvalues of  $(1-t) A_1 + t A_2$ have  opposite signs. It is not hard to construct such a matrix (see e.g~\cite{BLMZexample}, Example 1.3). Suppose $a_{12} = \beta t$ and $a_{21} = \beta (1-t)$ with $\beta > 0,$ so that $p_1 = (1-t), p_2 = t.$ Then, by Corollary \ref{cor:trace}, the Lyapunov exponents, $\lambda_1, \lambda_2$ (counted with their multiplicity) satisfy
 $$\lambda_1 + \lambda_2 = (1-t) \Tr(A^1) + t \Tr(A^2) < 0,$$
 while, it follows from Theorem 1.6 of \cite{BLMZexample}, that
 $\Lambda^+ = \Lambda^- > 0$ for $\beta$ sufficiently large. Hence (for large $\beta$)
 $$\lambda_2 < 0 < \lambda_1 = \Lambda^- = \Lambda^+.$$} \erem

\subsection{Uniqueness of average growth rate}
\label{sec:uniq}
In this section we discuss general conditions ensuring that  $$\Lambda^- = \Lambda^+ = \lambda_1.$$
A sufficient condition is given by {\em unique ergodicity} of $(\Theta, J),$ meaning that  ${\cal P}_{inv}^{(\Theta, J)}$ has cardinal one.
However, whenever $C_{M}$ is symmetric (i.e~ $C_M = -C_M$), for each   $\mu \in {\cal P}_{inv}^{(\Theta, J)}$ there is another (possibly equal) invariant measure $\mu^-$ given as the image measure of $\mu$ by the map $x,i \mapsto -x,i.$
 Indeed, it is easy to see that  $$[\mu P_t^{\Theta, J}]^- = \mu^- P_t^{\Theta, J}$$ for all $\mu \in {\cal P}(S^{d-1} \cap C_M \times E).$ This follows from  the  equivariance property $$G^i(-x) = -G^i(x)$$ satisfied by the $G^i$ (see equation \ref{eq:defGi}).
  Clearly $\Lambda(\mu) = \Lambda(\mu^-).$ Thus, when $C_M$ is symmetric, a (weaker than unique ergodicity) sufficient condition is that  the quotient space ${\cal P}_{erg}^{(\Theta, J)} /  \sim$ obtained by identification of $\mu$ with $\mu^-$ has cardinal one.

\bex[One dimensional systems]
{\rm
Suppose $d=1$ and $C_M = \RR.$ Thus  $S^{d-1} \cap C_M = \{ \pm 1\}$ and  ${\cal P}_{erg}^{(\Theta, J)} = \{\mu, \mu^-\}$ where $\mu^i(1) =\mu^{-,i}(-1) = p_i$ and $\mu^i(-1)  =\mu^{-,i}(1) =  0.$ Hence $ \Lambda^- = \Lambda^+ = \lambda_1 = \sum_i p_i a^i$ where $a^i =(F^i)'(0)$.
}
\eex

The two following results complement the previous discussion with practical conditions.

Set $\mathrm{G} = \{G^i\}_{i \in E}, \mathrm{G}_0 = \mathrm{G}, \mathrm{G}_{k+1} = \mathrm{G}_k \cup \{[G^i, V], V \in \mathrm{G}_k\}$ where $[,]$ is the Lie bracket operation.
Following \cite{BMZIHP}, we say that the {\em weak bracket}  condition holds at $p \in S^{d-1}$ provided the vector space spanned by the vectors $\{V(p) \: : V \in \cup_{k \geq 0}  \mathrm{G}_{k}\}$  has full rank (i.e $d-1$).

\bprop Assume there exists $p \in S^{d-1} \cap C_M$ such that
\label{th:uniqueLambda}
\bdes
\iti  The weak bracket condition holds at $p;$
\itii Either $p$ is $G$-accessible from $S^{d-1} \cap C_M$ or, $C_M$ is symmetric and $\{-p,p\}$ is $G$-accessible from $S^{d-1} \cap C_M.$
\edes
Then ${\cal P}_{inv}^{(\Theta, J)}$ in the first case, and ${\cal P}_{erg}^{(\Theta, J)} /  \sim$ in the second, has cardinal one. In particular
$$\Lambda^- = \Lambda^+ = \lambda_1.$$
\eprop
\prf Existence of an invariant probability follows from compactness and Feller continuity. By Theorem 1 in \cite{bakhtin&hurt} or Theorem 4.4 in \cite{BMZIHP} Condition $(i),$ and accessibility of $p$ imply that such a measure is unique (and absolutely continuous with respect to $dx \otimes  \sum_i \delta_i$). In case $C_M$ is symmetric and $\{-p,p\}$ accessible, let $S^{d-1} \cap C_M / \sim$ be the projective space obtained by identifying each point $x$ with the antipodal point $-x$ and $\pi : S^{d-1} \cap C_M  \mapsto S^{d-1} \cap C_M /\sim$ the quotient map. The PDMP $(\Theta, J)$ induces a PDMP $(\pi{\Theta},J) = (\pi(\Theta_t), J_t)$ on $S^{d-1} \cap C_M /\sim \times E$ for which $\pi(p)$ is accessible and at which the weak bracket condition holds. The preceding results applies again.
\qed
\bex[Two dimensional systems]
{\rm Suppose $d = 2, C_M = \RR^2$ and that one of the two following conditions is verified :
\bdes
\ita At least one matrix, say $A^1$,   has no real eigenvalues; or
\itb at least two matrices, say $A^1, A^2$  have no (nonzero) common eigenvector.
\edes
Then the assumptions, hence the conclusions, of Proposition \ref{th:uniqueLambda} hold.

Indeed, under condition $(a)$, the flow induced by $G^1$ is periodic on $S^1$ so that every point $p \in S^1$ satisfies the assumptions of Proposition \ref{th:uniqueLambda}. Under condition $(b)$, let $\alpha \leq \beta$ be the eigenvalues of $G^1$ and $u,v \in S^1$ corresponding eigenvectors. If $\alpha < \beta$ $\{v,-v\}$ is an attractor for the flow induced by $G^1$ whose basin is $S^1 \setminus \{u,-u\}.$ Since $G^2(u) \neq 0,$  $\{-v,v\}$ is $\{G^1,G^2\}$ accessible and since $G^2(v) \neq 0$ assumption $(i)$ of Proposition \ref{th:uniqueLambda} is satisfied at point $v.$ If $\alpha = \beta$ every trajectory of  the flow induced by $G^1$ converges either to $v$ or $-v$ and the preceding reasoning still applies.  }
\eex

The next proposition will be  useful in Section \ref{sec:epidemic}  for analyzing   random switching between {\em cooperative vector fields} and certain epidemiological models.
In case the matrices $A^i$ are irreducible, this proposition follows from the  Random Perron-Frobenius theorem as proved by Arnold, Demetrius and Gundlach in \cite{ADG}. However, to handle the weaker assumption $(iii),$ the proof needs to be adapted, but relies on the same ideas. Details are given in Section \ref{sec:appendix}. Recall (see remark \ref{rem:Mierc}) that a \textit{Metzler} matrix is a matrix  with nonnegative off-diagonal entries. We say that such a matrix is \textit{irreducible} if adding a sufficiently large multiple of the identity, the obtained matrix is a non-negative irreducible matrix in the usual sense. 
\bprop
\label{th:Perron} Assume that
\bdes
\iti $C_M = \RR^d_+,$
\itii For each $i \in E,$ $A^i$ is Metzler,
\itiii There exists $\alpha \in {\cal P}(E)$ (i.e $\alpha_i \geq 0, \sum_{i \in E} \alpha_i = 1$) such that
$$\overline{A} = \sum_{i \in E} \alpha_i A^i$$
 is  irreducible.
\edes
Then  ${\cal P}_{inv}^{(\Theta, J)}$ has cardinal one. In particular
$$\Lambda^- = \Lambda^+ = \lambda_1.$$
\eprop

 \subsection{Average growth rate under frequent switching}
The definition of  average growth rates (see equations (\ref{deflambdamu}) and (\ref{deflambdamumax})) involve the invariant measures of $(\Theta,J)$ whose explicit computation may prove highly difficult if not impossible.   However,
when switchings occur frequently, such measures can, by a standard averaging procedure, be estimated  by the invariant measures  of  the mean vector field; i.e the vector field  obtained by averaging.

More precisely, we have the following Lemma :

\blem
\label{lem:average1}
  Assume the switching rates are constant and depend on a small parameter $\varepsilon: a_{i,j}^{\varepsilon} = a_{i,j}/{\varepsilon} $ where $(a_{i,j})$ is an irreducible matrix with invariant probability $p$. Denote by $(\Theta^{\varepsilon},J^{\varepsilon})$ the associated PDMP
   given by~\eqref{dThetarho}, and for any $\varepsilon >0$, let $\mu^{\varepsilon}$ be an element of ${\cal P}_{inv}^{(\Theta^{\varepsilon}, J^{\varepsilon})}.$ Then, every limit point of $(\mu^{\varepsilon})_{\varepsilon > 0},$ in the limit $\eps \rar 0,$ is  of the form $\nu \otimes p$, where $\nu$ is an invariant probability measure of the flow induced by $G^p := \sum_i p_i G^i$.
  \label{lem:moyenne}
  \elem
The proof of this lemma follows from standard averaging results. Details are given in Section \ref{sec:appendix}.
An immediate corollary is :
\bcor
With the  hypotheses of Lemma \ref{lem:average1}, assume that the flow induced by $G^p $ admits a unique invariant measure $\nu$ on $S^{d-1} \cap C_M$. Denote by $\Lambda_{\varepsilon}^+$ and $\Lambda_{\varepsilon}^-$ the extremal growth rates of  $(\Theta^{\varepsilon},J^{\varepsilon})$. Then
$$ \lim_{ \varepsilon \to 0} \Lambda_{\varepsilon}^+ =  \lim_{ \varepsilon \to 0} \Lambda_{\varepsilon}^- = \sum_{i \in E} p_i \int_{S^{d-1} \cap C_M} \langle A^i \theta, \theta \rangle \nu(d\theta).$$
In particular, if $A^p := \sum_i p_i A^i$ is Metzler and irreducible, then it admits a unique eigenvector $\theta^p$ on $S^{d-1} \cap \RR^d_+$ and
$$ \lim_{ \varepsilon \to 0} \Lambda_{\varepsilon}^+ =  \lim_{ \varepsilon \to 0} \Lambda_{\varepsilon}^- = \langle A^p \theta^p, \theta^p \rangle = \lambda_{\max}(A^p). $$
\label{cor:moyenne}
\ecor
\section{The non linear system : Main results}
\label{sec:main}
\subsection{Extinction}
The first result is an {\em extinction} result.
\bthm
\label{th:extinct}
Assume $\Lambda^+ < 0.$ Let $0 < \alpha < - \Lambda^+.$ Then there exists a neighborhood ${\cal U}$ of $0$  and $\eta > 0$ such that for all $x \in {\cal U}$ and $i \in E$
$$\mathbb{P}^Z_{x,i}( \limsup_{t \rar \infty} \frac{1}{t} \log(\|X_t\|) \leq - \alpha) \geq \eta.$$
If furthermore $0$ is $\mathrm{F}$-accessible from $M,$ then for all $x \in M$ and $i \in E$
$$\mathbb{P}^Z_{x,i}( \limsup_{t \rar \infty} \frac{1}{t} \log(\|X_t\|) \leq \Lambda^+) = 1.$$
\ethm
\subsection{Persistence}
The next results are {\em persistence} results obtained under the assumption that $\Lambda^- > 0.$

We let $$\Pi_t = \frac{1}{t} \int_0^t \delta_{Z_s} ds \in {\cal P}(M \times E)$$ denote  the {\em empirical occupation measure} of the process $Z.$ For every Borel set $A \subset M \times E$
$$\Pi_t(A) = \frac{1}{t} \int_0^t \Ind_{\{Z_s \in A\}} ds$$ is then the proportion of the time spent by $Z$ in $A$ up to time $t.$

We let $M^* = M \setminus \{0\}.$
\bthm

\label{th:persist1}
Assume $\Lambda^- > 0.$ Then the following assertions hold:
\bdes \iti For all $\eps > 0$ there exists $r > 0$ such that for all $x \in M^*$, $i \in E$, $\mathbb{P}^Z_{x,i}$ almost surely,
$$\limsup_{t \rar \infty} \Pi_t(B(0,r) \times E) \leq \eps.$$
In particular, for all $x \in M^*,$ $\mathbb{P}^Z_{x,i}$ almost surely,  every limit point (for the weak* topology) of $(\Pi_t)$ belongs to ${\cal P}_{inv}^Z \cap {\cal P}( M^* \times E).$
\itii There exist positive constants $\theta, K$ such that for all $\mu \in {\cal P}_{inv}^Z \cap {\cal P}( M^* \times E)$
$$\sum_{i \in E} \int \|x\|^{-\theta} \mu^i(dx) \leq K.$$
\itiii Let $\eps > 0$ and $\tau^{\eps}$ be the stopping time defined by
$$\tau^{\eps} = \inf \{t \geq 0 : \: \|X_t\| \geq \eps \}.$$ There exist $\eps > 0$,  $b > 1$ and $c > 0$ such that
for all $x \in M^*$ and $i \in E$,
$$\mathbb{E}_{x,i}^Z(b^{\tau^{\eps}}) \leq c (1 + \|x\|^{-\theta}).$$
\edes
\ethm
 Set $\mathrm{F}_0 = \mathrm{F} = \{F^i\}_{i \in E}$ and $\mathrm{F}_{k+1} = \mathrm{F}_k \cup \{[F^i, V], V \in \mathrm{F}_k\}$ where $[, ]$ is the Lie bracket operation.
 We say (compare to Section \ref{sec:uniq}) that  the {\em weak bracket} condition holds at $p \in M$ provided the vector space spanned by the vectors $\{V(p) \: : V \in \cup_{k \geq 0}  \mathrm{F}_{k}\}$  has full rank.
 We let $\mathbf{Leb}$ denote the Lebesgue measure on $\RR^d.$
\bthm
 \label{th:persist2} In addition to the assumption $\Lambda^- > 0,$  assume that there exists a point $p \in M^*$ $\mathrm{F}$-accessible from $M^*$ at which the weak bracket condition holds. Then
 \bdes
 \iti
The set ${\cal P}_{inv}^Z \cap {\cal P}( M^* \times E)$ reduces to a single element, denoted $\Pi$;
\itii  $\Pi$ is absolutely continuous with respect to $\mathbf{Leb} \otimes (\sum_{i \in E} \delta_i)$;
\itiii For all $x \in M^*$ and $i \in E$, $$\lim_{t \rar \infty} \Pi_t = \Pi$$ $\mathbb{P}^Z_{x,i}$ almost surely. \edes
\ethm
In order to get a convergence in distribution of the process $(Z_t)_{t \geq 0}$, the weak bracket condition needs to be strengthened.
Set  $\mathcal{F}_{0} = \{F^i  - F^j\: : i, j = 1, \ldots m\}$ and $\mathcal{F}_{k+1}  =  \mathcal{F}_k \cup \{[F^i , V ] \: : V \in \mathcal{F}_k\}.$
We say that the {\em strong bracket} condition holds at $p \in M$ provided the vector space spanned by the vectors $\{V(p) \: : V \in \cup_{k \geq 0}  \mathcal{F}_{k}\}$  has full rank.

Given $\mu, \nu \in {\cal P}(M \times E),$ the {\em total variation distance} between $\mu$ and $\nu$ is defined as
$$\|\mu - \nu\|_{TV} = \sup |\mu(A) - \nu(A)|$$ where the supremum is taken over all Borel sets $A \subset M \times E.$
\bthm
 \label{th:persist3}
Under the conditions of the preceding theorem, assume furthermore that one the two following holds :
\bdes
\iti The weak bracket condition is strengthened to the strong bracket condition; or
\itii There exist $\alpha_1,\ldots,\alpha_N \in \RR$ with $\sum \alpha_i = 1$ and a point $e^{\star} \in M^*$ $\mathrm{F}$-accessible from $M^*$ such that  $\sum \alpha_i F^i(e^{\star}) = 0$.
\edes
  Then there exist $\kappa, \theta > 0$ such that for all $x \in M^*$ and $i \in E$,
$$\|\mathbb{P}^Z_{x,i}(Z_t \in \cdot ) - \Pi\|_{TV} = \| \delta_{x,i} P_t^Z - \Pi \|_{TV} \leq const. (1 + \|x\|^{-\theta}) e^{-\kappa t}.$$
\ethm

\subsection{The noncompact case}
\label{sec:noncompact}
We briefly discus here the situation where $M$ is not compact. First,  note that all the results given in section \ref{sec:linear} still hold, because they only deal with the linearised system. Next, local statements remain true without additional assumption by a localisation argument. Namely  :

\bthm \
\label{th:noncompactlocal}
\begin{enumerate}
\item Assume $\Lambda^+ < 0.$ Let $0 < \alpha < - \Lambda^+.$ Then there exists a neighborhood ${\cal U}$ of $0$  and $\eta > 0$ such that for all $x \in {\cal U}$ and $i \in E$
$$\mathbb{P}^Z_{x,i}( \limsup_{t \rar \infty} \frac{1}{t} \log(\|X_t\|) \leq - \alpha) \geq \eta.$$
\item Assume $\Lambda^- > 0.$ Then there exist $\eps > 0$,  $b > 1$ and $c > 0$ such that
for all $x \in M^*$ and $i \in E$,
$$\mathbb{E}_{x,i}^Z(b^{\tau^{\eps}}) \leq c (1 + \|x\|^{-\theta}).$$
\end{enumerate}
\ethm 
To extend the global results stated above  , we make the additional assumption that  the jumps rates are bounded and  that there exists a Lyapunov function, controlling the behaviour of the process at infinity. 

\begin{hypothesis}
\label{boundedjump}
The jumps rate are bounded : $$ \sup_{x \in M} \max_{i,j} a_{ij}(x) < \infty.$$
\end{hypothesis}
For a function $f : M \times E \to \RR$, we denote by $\Gamma f$ the function defined by :
$$ \Gamma f (x,i) = \sum_{j \in E} a_{ij}(x) \left( f(x,j) - f(x,i) \right)^2.$$
We also let $C^1_c$ denote the space of functions $f : M \times E \to \RR$ that are constant outside a compact set and $C^1$ in the first variable.
\begin{hypothesis}
\label{Lyapunovinfinity}
There exists a continuous function $W : M \times E \to \RR_+$ with $\lim_{\|x\| \to \infty} W(x,i)=\infty$, a continuous function $LW : M \times E \to \RR_+$, $\alpha > 0$ and $C \geq 0$ such that 
\begin{enumerate}
\item[\textbf{(i)}] For every compact set $K \subset M$, there exists $W_K \in C_c^1$ such that
\begin{enumerate}
\item[\textbf{(a)}] $W|_K = W_K|_K$ and ${\cal L}W_K|_K = LW|_K$,
\item[\textbf{(b)}] For all $x \in M$, $\sup \{ P_t( \Gamma W_K), \quad t \geq 0, \quad K \quad \text{compact} \} < \infty$
\end{enumerate}
\item[\textbf{(ii)}] $$ L W \leq - \alpha W + C. $$
\end{enumerate}

\end{hypothesis}

\bthm
\label{th:noncompactglobal}
Under Hypotheses \ref{Lyapunovinfinity} and \ref{boundedjump}, Theorems \ref{th:extinct}, \ref{th:persist1} and \ref{th:persist2} are still valid. Moreover, Theorem \ref{th:persist3} is true, but with the following estimate :
$$ \| \delta_{x,i} P_t^Z - \Pi \|_{TV} \leq const. (1 + W(x) + \|x\|^{-\theta} ) e^{-\kappa t}.$$
\ethm

\bex {\rm
We consider a random switching between two linear systems given by $ 2 \times 2$ Metzler matrices $A^0$ and $A^1$, with transition rate $a_{i,1-i}(x)$. We assume that $A^0$ has two distinct positive eigenvalues $\lambda_1 > \lambda_2$ and is irreducible, whereas $A^1$ is of the form $$A^1 =\begin{pmatrix}
- c & 0\\
0 & -d
\end{pmatrix},$$
with $0 < c < d$.
Since the eigenvalues of $A^0$ are positive, there is no invariant compact set for $\Psi^0$, nor for the PDMP. Moreover, $A^0$ and $A^1$ being Metzler,  $M=\RR^2_{+}$ is positively invariant for $(X_t)_{t \geq 0}$. If the jump rates were constant in $x$, the process would either converge to $0$ or to infinity. To ensure positive recurrence on $M^*$,  we assume that the transition rates are such that, near the origin, $I_t$ spends more time in state $0$ :

\beq
\label{eqnc1}
a_{10}(0) -  \frac{d}{\lambda_2} a_{01}(0)>0;
\eeq
While near infinity, it spends more time in state $1$ :

\beq
\label{eqnc2}
 \limsup_{ \|x\| \to \infty} \left( a_{10}(x) -  \frac{c}{\lambda_1} a_{01}(x) \right) < 0.
\eeq
More precisely, we have the following :

\bprop
Assume that the jumps rates are bounded and that conditions \eqref{eqnc1} and \eqref{eqnc2} hold. Then there exists a unique invariant probability $\Pi \in \mathcal{P}(M^* \times E)$ and there exists  $\kappa, \theta, q > 0$ such that for all $x \in M^*$ and $i \in E$,
$$\|\mathbb{P}^Z_{x,i}(Z_t \in \cdot ) - \Pi\|_{TV}  \leq const. (1 + \|x\|^q + \|x\|^{-\theta}) e^{-\kappa t}.$$
\eprop  

\prf By Theorem \ref{th:Perron}, $\Lambda^+ = \Lambda^-:=\Lambda$, and by Corollary \ref{cor:trace},  $$\Lambda \geq \frac{1}{2}(p_0 \Tr(A^0) + p_1 \Tr(A^1) \geq \lambda_2 p_0 - d p_1.$$ 
Moreover, it is easy to check that $p_0 = \frac{a_{10}(0)}{a_{10}(0)+a_{01}(0)}$ and $p_1 = \frac{a_{01}(0)}{a_{10}(0)+a_{01}(0)}$. Hence, if $a_{10}(0) > \frac{d}{\lambda_2} a_{01}(0)$,  then $\Lambda > 0$. Now we show that we can construct a Lyapunov function at infinity. Let $q > 0$ and $\beta_0, \beta_1 > 0$ and  define, for all $(x,i) \in M \times E$, $W_q(x,i) = \beta_i \| x \|^q$. Formally, we have $$ { \cal L}W_q(x,i) = q \beta_i \langle A_i x, x \rangle \| x \| ^{p-2} + a_{i, 1-i}(x) ( \beta_{1-i} - \beta_i) \|x \|^q.$$
 By assumption on $A^0$ and $A^1$, $\langle A_0 x, x \rangle \leq \lambda_1 \| x \|^2$ and $\langle A_1 x, x \rangle \leq - c \| x \|^2
$. Hence,  
$$ { \cal L}W_q(x,i) \leq  \left(- \alpha(i) q \beta_i  + a_{i, 1-i}(x) ( \beta_{1-i} - \beta_i) \right) \|x \|^q,$$
where $\alpha(0) = - \lambda_1$ and $\alpha(1) = c$. First we prove that we can choose $\beta_0$ and $\beta_1$ such that $W_q$ satisfies point \textbf{(ii)} of Hypothesis \ref{Lyapunovinfinity} for all $q$ small enough. Then we prove that we can choose $q$ such that point \textbf{(i-b)} holds. By assumption \eqref{eqnc2}, there exists $\eps>0$ and $K>0$ such that, for all $x \in M$ with $\|x\| \geq K$,  $a_{10}(x) \leq  \frac{c}{\lambda_1} a_{01}(x) - \eps$. This implies that, for $q$ small enough, there exists $\alpha_q$ such that $a_{10}(x)( \frac{\alpha_q}{\lambda_1}+q) -  (\frac{c}{\lambda_1} -  \frac{\alpha_q}{\lambda_1}) a_{01}(x) - q \alpha_q+ cq^2 \leq 0$, which yields $$\sup_{ \|x\| \geq K} \frac{a_{01}(x)+\alpha_q}{a_{01}(x)-\lambda_1 q} \leq \inf_{ \|x\| \geq K} \frac{-\alpha_q+cq}{a_{10}(x)}+1.$$ Now we choose $\beta_1=1$ and $\beta_0$ such that $$\sup_{ \|x\| \geq K} \frac{a_{01}(x)+\alpha_p}{a_{01}(x)-\lambda_1 q} \leq \beta_0 \leq \inf_{ \|x\| \geq K} \frac{-\alpha_q+cq}{a_{10}(x)}+1.$$ Thus, for $\|x\| \geq K$, $ - \alpha(i) q \beta_i  + a_{i, 1-i}(x) ( \beta_{1-i} - \beta_i) \leq - \alpha_q$. In particular, for all for $\|x\| \geq K$, ${ \cal L}W_q(x,i) \leq - \alpha_q W_q(x,i)$. Since ${ \cal L}W_q$ is bounded for $\|x\| \leq K$, then ${ \cal L}W_q \leq - \alpha_q W_q + C$ for some constant $C > 0$ (depending on $q>0$). This has the consequence (see \cite[Theorem 2.1]{Ben14}) that for all $t \geq 0$, 
\beq
 P_t W_q \leq e^{- \alpha_q t} \left( W_q - \frac{C}{\alpha_q}\right) + \frac{C}{\alpha_q}.
 \label{eq:ptlyap}
 \eeq 
The computation of $\Gamma$ gives $$\Gamma W_q(x,i) = a_{i,1-i}(x) ( \beta_0 - \beta_1)^2 \| x \|^{2 q},$$ hence $$ \Gamma W_q \leq \tilde{C}_q W_{2q}$$ for some constant $\tilde{C}_q > 0$. Hence, choosing $p$ small enough so that \eqref{eq:ptlyap} holds for $2q$, one has $$ \sup_{ t \geq 0} P_t \left( \Gamma W_q \right) \leq \tilde{C}_q \sup_{ t \geq 0} P_t W_{2q} \leq W_{2q},$$ which proves \textbf{(i-b)}. It remains to show that there exist accessible points at which the strong bracket condition holds. Set $F^0(x) = A^0x$ and $F^1(x) = A^1x$ the vector fields associated to $A^0$ and $A^1$. There exist $\alpha, \beta, \gamma, \delta$, with $\beta, \gamma > 0$ such that $F^0(x,y) = ( \alpha x + \beta y, \gamma x + \delta y)$. Straightforward computations show that $$\det( F^0-F^1, [F^0,F^1])(x,y) = (d-c)(2 \beta \gamma x y + \beta ( d + \delta) y^2 + \gamma ( \alpha + c) x^2).$$ Since $\beta, \gamma > 0$, this polynomial is non identically null. To conclude, we prove that there exists an open set of accessible points. Let $v \in \RR^2_{++}$ be the Perron eigenvector associated with $A^0$. We claim that $\RR_+v$ and therefore $\gamma_1^+(\RR_+ v) = \overline{\cup_{ t \geq 0} \Psi^1_t(\RR_+ v)} $ are accessible. One can check that for all $y \in \RR_+ v$ and all $\eps > 0$, there exists $\eta > 0$ such that for all $x \in M^*$ with $ \|x \| < \eta$, there exists $t \geq 0$ such that $\| \Psi^0_t(x) - y \| < \eps$. Since $0$ is accessible following $F^1$, this makes $y$ accessible. Hence, $\gamma_1^+(\RR_+ v)$ is accessible and Theorem \ref{th:noncompactglobal} applies. \qed

}
\eex

\section{Epidemic Models in Fluctuating Environment}
\label{sec:epidemic}
We discuss here some implications of our results to certain epidemics models evolving in a randomly fluctuating environment.

Forty years ago,  Lajmanovich and Yorke  in a influential   paper \cite{LajYorke},   proposed and analyzed a  deterministic SIS (susceptible-infectious-susceptible) model of infection,  describing the evolution of a disease that does not confer immunity, in a population structured in $d$ groups.
The model is given by a differential equation on $[0,1]^d$ (the unit cube of $\RR^d$)  having  the form
\beq
\label{eq:LY}
\frac{dx_i}{dt} = (1-x_i) (\sum_{j = 1}^d C_{ij} x_j) - D_i x_i\, ,   \: i = 1, \ldots d,
\eeq
where  $C = (C_{ij})$ is an irreducible matrix with nonnegative entries and $D_i > 0.$
Here $0 \leq x_i \leq 1$ represents the proportion of infected  individuals in group $i;$ $D_i$ is the intrinsic cure rate in group $i$ and $C_{ij} \geq 0$ is the rate at which group $i$ transmits the infection to group $j.$ Irreducibility of $C$ implies that each group indirectly affects the other groups. By a classical mean field approximation procedure, (\ref{eq:LY}) can be derived from a finite population  model, in the limit of an infinite population (see Bena{\"i}m and Hirsch \cite{BH99}).

Here and throughout, for any matrix $A$ we let  $\lambda(A)$ denote the largest real part of the eigenvalues of  $A.$ A matrix $A$ is called \textit{Hurwitz} provided $\lambda(A) < 0$.
Lajmanovich and Yorke \cite{LajYorke} prove the following result:
\bthm[Lajmanovich and Yorke, \cite{LajYorke}]
 \label{th:LajYorke} Let  $A = C-\mathsf{diag}(D).$

 If $\lambda(A) \leq 0,$  $0$ is globally asymptotically stable for the semiflow induced by (\ref{eq:LY}) on $[0,1]^d.$

If $\lambda(A) > 0$ there exists another  equilibrium $x^* \in ]0,1[^d$ whose basin of attraction is $[0,1]^d \setminus \{0\}.$
\ethm

In this epidemiological framework, $0$ is called the \textit{disease free equilibrium}, and the point $x^*$, when it exists, the \textit{endemic equilibrium}.
It turns out that such a  dichotomic behavior  is very robust to the perturbations of the model and can be obtained under a very general set of assumptions, using Hirsch's  theory of {\em cooperative differential equations}.

We let $\RR^d_{++}$ denote the  interior of the non negative orthant $\RR^d_+.$
For $x, y \in \RR^d$ we write $x \leq y$ (or $y \geq x$) if $y - x \in \RR^d_+; x < y$ if $x \leq y$ and $x \neq y;$ and $x << y$ if $y -x \in  \RR^d_{++}.$

Following \cite{BH99} (especially Section 3), we call a map $F : [0,1]^d \mapsto \RR^d$  an {\em epidemic vector field} if it is continuously differentiable\footnote{by this we mean that $F$ can be extended to a $C^1$ vector field on $\RR^d.$} and satisfies the following set of conditions:
\bdes
\item[E1] $F(0) = 0;$
\item[E2] $x_i = 1 \Rightarrow F_i(x) < 0;$
\item[E3] $F$ is {\em cooperative} i.e the Jacobian matrix  $DF(x)$ is Metzler for all $x \in [0,1]^d;$
\item[E4] $F$ is {\em irreducible} on   $[0,1)^d$ i.e~  $DF(x)$ is irreducible for all $x \in [0,1)^d;$
\item[E5] $F$ is strongly {\em sub-homogeneous} on $(0,1)^d$ i.e~ $F(\lambda x) << \lambda F(x)$ for all $\lambda > 1$ and $x \in (0,1)^d.$
\edes
It is easy to verify that the {\em Lajmanovich and Yorke vector field} (given by the right hand side of (\ref{eq:LY})) satisfies these conditions.

Let $\Psi = \{\Psi_t\}$ denote the local flow induced by $F.$ Condition $E3$ has the important consequence that for all $t \geq 0$  $\Psi_t$  is {\em monotone} for the partial ordering $\leq.$ That is  $\Psi_t(x) \leq \Psi_t(y)$ if $x \leq y.$  In particular, by $E1,$ $\Psi_t(x) \geq 0$ for all $x \geq 0.$ Combined with $E2$ this shows that
$[0,1]^d$ is positively invariant under $\Psi.$

The following result shows that trajectories of $\Psi$ behave exactly like the trajectories of the Lajmanovich and Yorke system.
The first assertion  was stated in (\cite{BH99}, Theorem 3.2) but its proof is a consequence of more general results due to Hirsch (in particular Theorems 3.1 and 5.5  in \cite{Hirsch94}).
\bthm
\label{th:LajYorkeGen}
Let $F$ be an epidemic vector field and $\Psi = \{\Psi_t\}_{t \geq 0}$ the induced semiflow on $[0,1]^d.$ Then
\bdes
\iti (Hirsch, \cite{Hirsch94}) Either $0$ is globally asymptotically stable for $\Psi$; or there exists another equilibrium $x^* \in ]0,1[^d$ whose basin of attraction is $[0,1]^d \setminus \{0\}.$
\itii Let $A = DF(0).$ Then $0$ is globally asymptotically stable if and only if $\lambda(A) \leq 0.$
\edes
\ethm
\prf As already mentioned, $(i)$ follows from \cite{Hirsch94}, Theorems 3.1 and 5.5. We detail the proof of $(ii)$. If $\lambda(A) < 0,$ then $0$ is linearly stable hence globally stable by $(i).$ If $\lambda(A) > 0,$ there exists, by irreducibility and Perron Frobenius theorem, $x_0 >> 0$ such that $Ax_0 = \lambda(A) x_0 >> 0.$ Hence $F(\eps x_0) >> 0$ for $\eps$ small enough, because $\frac{F(\eps x_0)}{\eps} \rar Ax_0$ as $\eps \rar 0.$ Consequently $\{x: \: x \geq \eps x_0\}$ is positively invariant and $0$ cannot be asymptotically stable.

It remains to show that $0$ is asymptotically stable when  $\lambda(A) = 0.$ Suppose the contrary. By $(i)$ there exists another equilibrium $x^* >>0.$ Set $y^* = x^*/2.$ By strong subhomogeneity, $0 = F(x^*) << 2 F(y^*).$ Let $F_{\eps}(x) = F(x) - \eps x.$ For all $\eps > 0,$  $F_{\eps}$ is an epidemic vector field and  $0$ is linearly stable for $F_{\eps}$ (because $\lambda (DF_{\eps}(0))  = - \eps$). On the other hand, for $\eps$ small enough, $0 << F_{\eps}(y^*)$ so that the set $\{y: \: y \geq y^*\}$ is  positively invariant by $F_{\eps}.$  A contradiction.
\qed
\subsection{Fluctuating environment}
We consider a PDMP $Z = (X,I)$ as defined in Section \ref{sec:intro}, under the assumptions that:
\bdes
\item[E'1] $M = [0,1]^d;$
\item[E'2] For all $i \in E,$ $A^i = DF^i(0)$ is Metzler;
\item[E'3] There exists $\alpha \in {\cal P}(E)$ such that  the convex combination
  $\overline{A} = \sum_{i \in E} \alpha_i A^i$ is irreducible.
\edes

Observe that these conditions are automatically satisfied if $\mathrm{F}  = \{F^i\}_{i \in E}$ consists of epidemic vector fields but are clearly much weaker.

Relying on Proposition \ref{th:Perron}, we let $\lambda_1 = \Lambda^+ = \Lambda^-$ denote the top Lyapunov exponent of the linearized system.
\bthm
  \label{th:EpidemicExtinct}
  Assume $\lambda_1 < 0$ and that one of the following two conditions holds:
  \bdes
  \ita The jump rates are constant (i.e~$a_{ij}(x) = a_{ij})$  and the $F^i$ are  epidemic; or
  \itb There exists $\beta \in {\cal P}(E)$ such that $\overline{F} = \sum_{i} \beta_i F^i$  is  epidemic and $$\lambda(\sum_i \beta_i A^i) \leq 0.$$
  \edes
  Then for all $x \in M^*$ and $i \in E$,
  $$\mathbb{P}^Z_{x,i}(\limsup \frac{\log(\|X_t\|)}{t} \leq \lambda_1) = 1.$$
  \ethm
  \prf We first prove the result under condition $(a).$  Recall (see Section \ref{sec:lyapou}) that $\Omega$ stands for  $D(\RR^+, E).$ For each $\omega \in \Omega$ and  $x \in [0,1]^d$ let $$t \mapsto \Psi(t,\omega)(x)$$ be the solution to the non autonomous differential equation $$\dot{y} = F^{\omega_t}(y),$$ with initial condition $y(0) = x.$    By conditions $E3$ and $E5$  each flow $\Psi^i$ is monotone and subhomogenous (see e.g~\cite{Hirsch94}, Theorem 3.1). The composition of monotone subhomogeneous mappings being monotone and subhomogeneous, $\Psi(t,\omega)$ is monotone and subhomogeneous for all $t \geq 0$ and $\omega \in \Omega.$
  Thus, for all $\eps > 0$ and $\|x\| > \eps$
  \beq
  \label{eq:Psisubhom}
  \Psi(t,\omega)(x) \leq \frac{\|x\|}{\eps} \Psi(t,\omega)(\frac{\eps}{\|x\|}x).
  \eeq
  Under the assumption that the jump rates are constant, $\mathbb{P}^Z_{x,i}$ is the image measure of $\mathbb{P}^J_{i}$ by the map
  $$\omega \mapsto (\omega, (\Psi(t,\omega)(x))_{t \geq 0}).$$ Therefore, by Theorem \ref{th:extinct}, there exists $\eta, \eps > 0$ such that for all $x \in B(0,\eps)$
  \beq
  \label{eq:couplingsubhom}
  \mathbb{P}^Z_{x,i} (\limsup_{t \rar \infty} \frac{\log(\|X_t\|)}{t} \leq \lambda_1) = \\
  \mathbb{P}^J_{i} (\limsup_{t \rar \infty} \frac{\log(\|\Psi(t,\omega)(x)\|)}{t} \leq \lambda_1) \geq \eta.
  \eeq
  Combined with (\ref{eq:Psisubhom}), this proves that (\ref{eq:couplingsubhom}) holds true not only for $x \in B(0,\eps)$ but for all
  $x \in [0,1]^d.$ A standard application of the Markov property then implies the result.

  Under condition $(b)$, it follows from Theorem \ref{th:LajYorkeGen}, that $0$ is $\mathrm{F}$-accessible from $M$, and the result follows from Theorem \ref{th:extinct}.
  \qed
  \brem{\rm The assumption made in case $(a)$ that the $F^i$ are epidemic  can be weakened. The proof shows that irreducibility of $F^i$ is unnecessary and that strong subhomogeneity can be weakened to subhomogeneity.}
  \erem
  \brem
  {\rm Case $(a)$ (and its proof) can be related with the  results obtained by
 Chueshov in \cite{chueshov}, for SIS models with random coefficients (see \cite[Section 5.7.2]{chueshov}) and, more generally,  for  monotone subhomogeneous  random dynamical systems.
  Note, however, that in comparison with Chueshov's approach,
in case $(b),$ there is no assumption that the $F^i$s are monotone nor subhomogeneous.}
  \erem
\bex[Fluctuations may promote cure] {\rm
We give here a simple example consisting of two  Lajmanovich-Yorke vector fields
 modeling the evolution of an endemic disease (each vector field possesses an endemic equilibrium) but such that a random switching between the dynamics leads to the extinction of the disease.

 Suppose $d = 2, E = \{0,1\}.$ Let  $F^0, F^1$ be the Lajmanovich-Yorke vector fields respectively given by  $$C^0 = \begin{pmatrix}
   2 & 1 \\
   1 & 1
\end{pmatrix}, \: D^0  = \begin{pmatrix}
   6 \\
    1
\end{pmatrix},$$ and $$C^1 = \begin{pmatrix}
   1 & 1 \\
   1 & 3
\end{pmatrix}, \:   D^1 = \begin{pmatrix}
   1 \\
    7
\end{pmatrix}.$$
One can easily  check that $$\lambda(A^0) = \lambda(A^1) = \sqrt{5}-2 >0,$$ so that for each $F^i,$ there is an endemic equilibrium and the  disease free equilibrium is  a repellor.
On the other hand, $$\lambda( \frac{A^0 + A^1}{2} )= - 1 < 0,$$ so that the disease free equilibrium is a global attractor of the average vector field $\overline{F} = \frac{1}{2}(F^1 + F^2).$ Consider now the PDMP  given by constant switching rates   $$a_{0,1} = a_{1,0} = \beta, a_{0,0}=a_{1,1}=0.$$   By Corollary~\ref{cor:moyenne}, this implies that $\lambda_1 < 0$ provided  $\beta$ is sufficiently large. Thus the conclusion of Theorem \ref{th:EpidemicExtinct} holds.

\begin{figure}
   \begin{minipage}[c]{.46\linewidth}
      \includegraphics[scale=0.3]{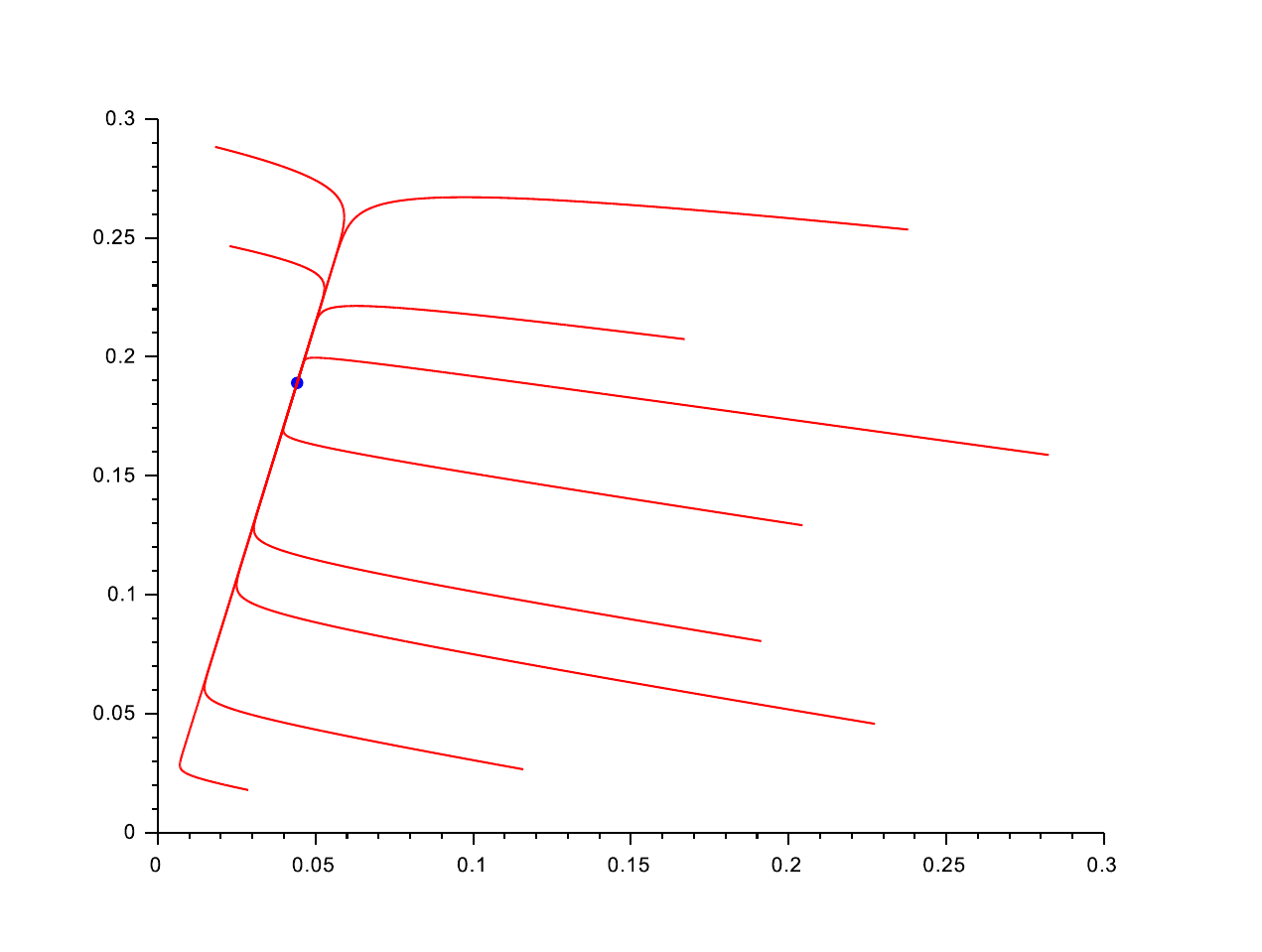}
   \end{minipage} \hfill
   \begin{minipage}[c]{.46\linewidth}
      \includegraphics[scale=0.3]{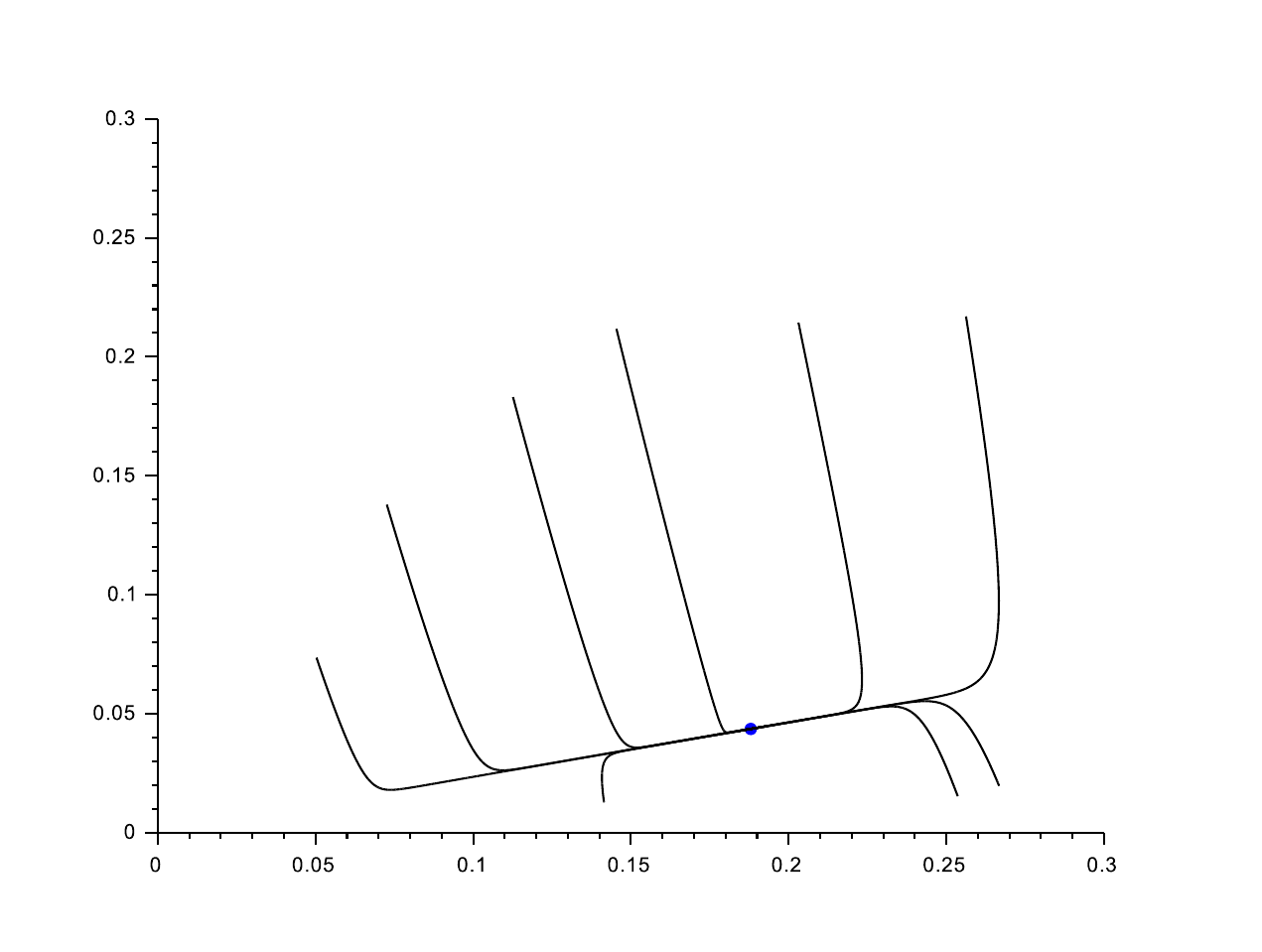}
   \end{minipage}
\caption{Example \ref{AinsCosta}, phase portrait of $F^0$ and $F^1$ }
\end{figure}

\begin{figure}[!h]
\centering
\includegraphics[scale=0.4]{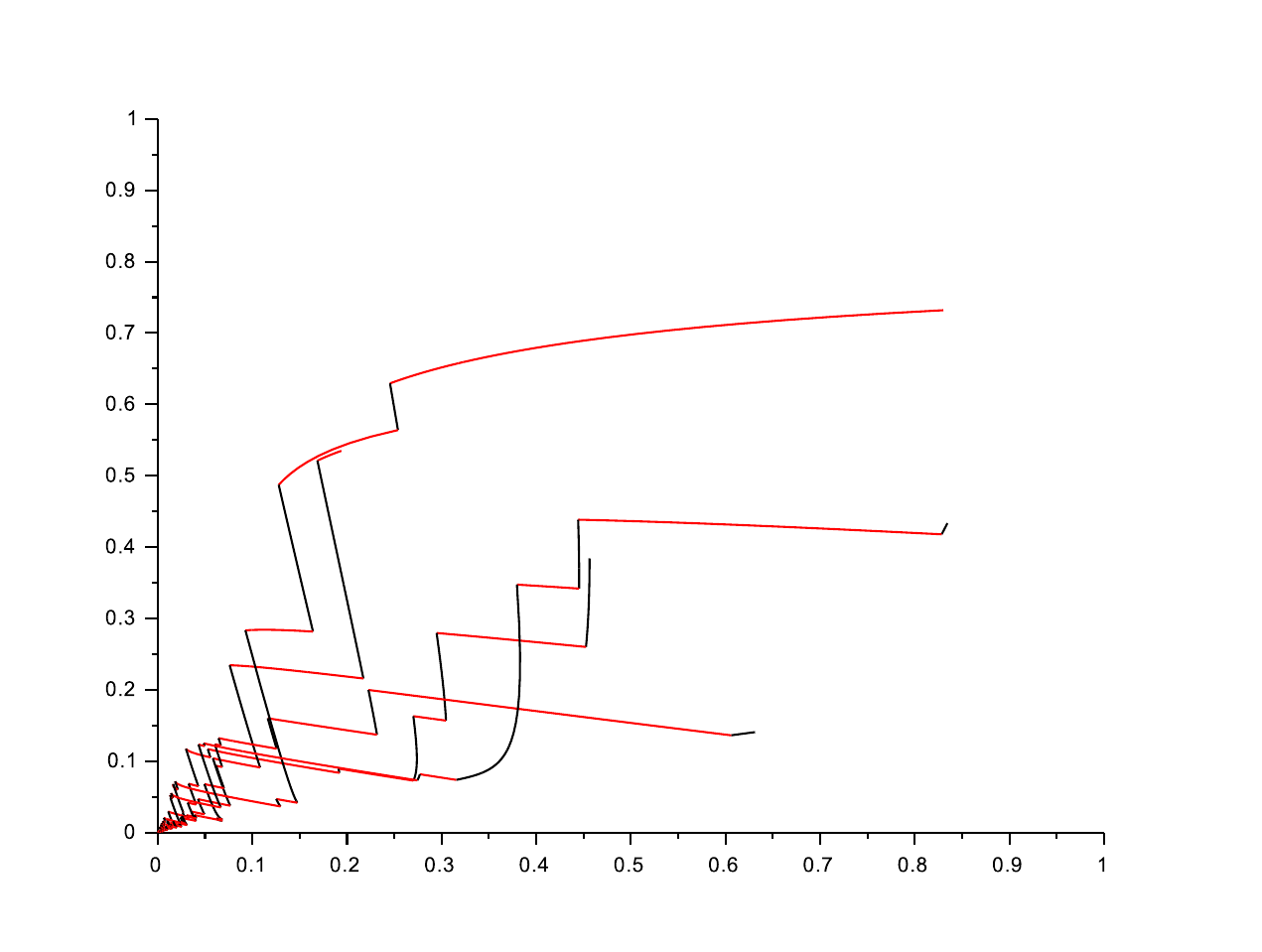}
\caption{Example \ref{AinsCosta}, some trajectories of $(X_t)$ for $\beta = 20$}
\end{figure}

  \label{AinsCosta}
}\eex

  \bex[Fluctuations may promote infection]{\rm We give here another simple example consisting of two Lajmanovich-Yorke vector fields for which the disease dies out,
  but such that  a random switching between the dynamics leads to the persistence of the disease.

  With the notation of Example~\ref{AinsCosta}, assume now that
 $$C^0 = \begin{pmatrix}
   1 & 4 \\
   \frac{1}{16} & 1
\end{pmatrix}, \; D^0  = \begin{pmatrix}
   2  \\
    2
\end{pmatrix},$$ and $$C^1 = \begin{pmatrix}
   2 &  \frac{1}{16} \\
   4 & 2
\end{pmatrix},   \:  D^1 = \begin{pmatrix}
   3  \\
   3
\end{pmatrix}.$$
  Straightforward computation shows that $$\lambda(A^0) = \lambda(A^1) = -1/2 < 0,$$  $$\lambda(\frac{A^0 + A^1}{2}) = 33/32 > 0,$$ and that the endemic equilibrium of $\overline{F}$  is the point $x^{\star} = (33/113, 33/113).$ Then $x^{\star}$ is $F$ - accessible and one can easily check that the strong bracket condition holds at $x^{\star}$. Thus, for $\beta$ sufficiently large, this implies by Corollary~\ref{cor:moyenne} and Theorem \ref{th:persist3} the exponential convergence in total variation of the distribution of $Z_t$ (whenever $X_0 \neq 0$) towards a unique distribution $\Pi$ absolutely continuous with respect to $\mathsf{Leb} \otimes \sum_{i \in E} \delta_i$ and satisfying the tail condition given by Theorem \ref{th:persist1} (ii). Furthermore, it follows from (\cite{BMZIHP}, Proposition 3.1) that the topological support of $\Pi$ writes
  $\Gamma \times E$ where $\Gamma$ is a compact connected set containing both $0$ and $x^{\star}$, and whose interior is dense in $\Gamma$.

\begin{figure}
   \begin{minipage}[c]{.46\linewidth}
      \includegraphics[scale=0.3]{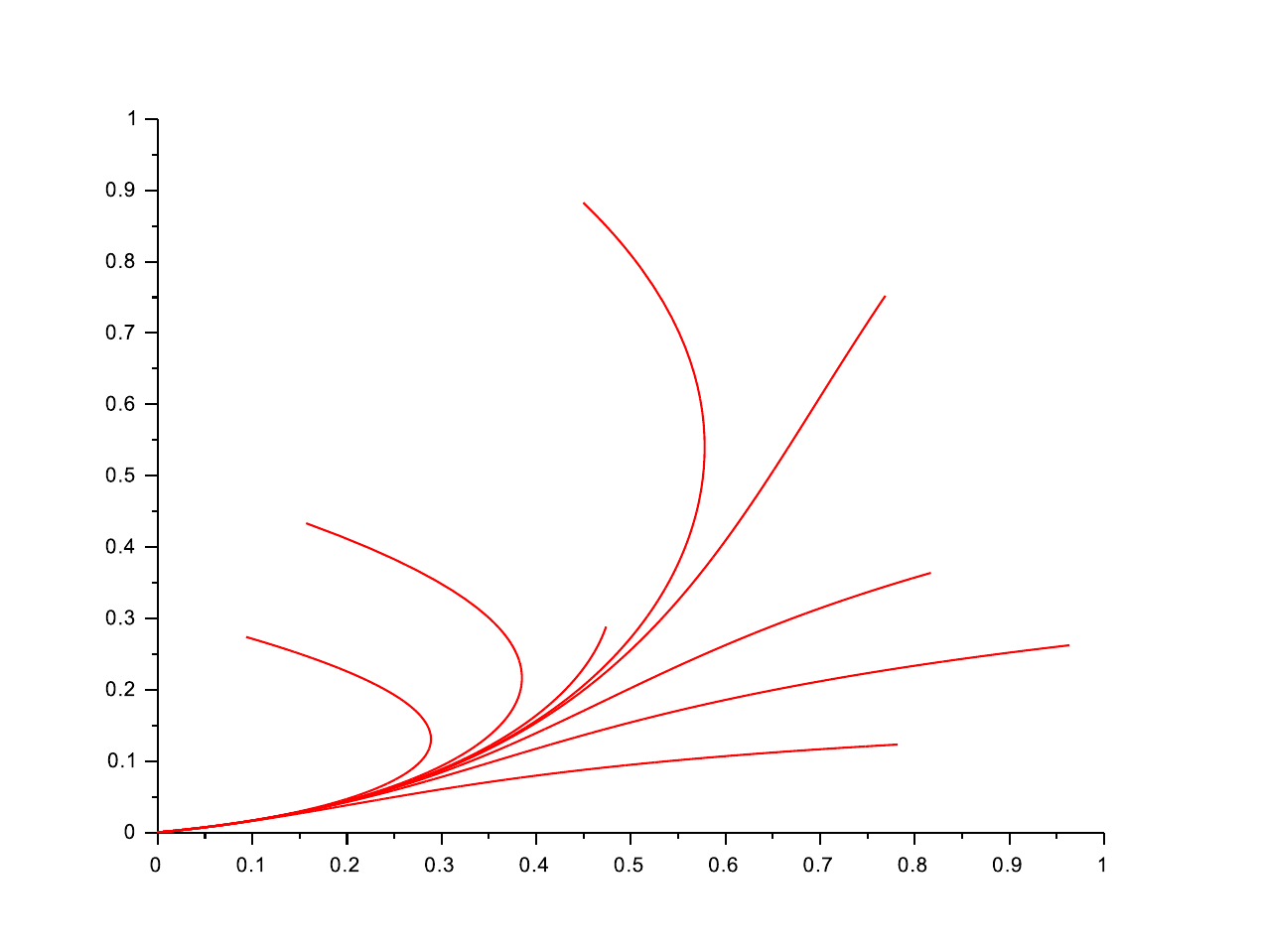}
   \end{minipage} \hfill
   \begin{minipage}[c]{.46\linewidth}
      \includegraphics[scale=0.3]{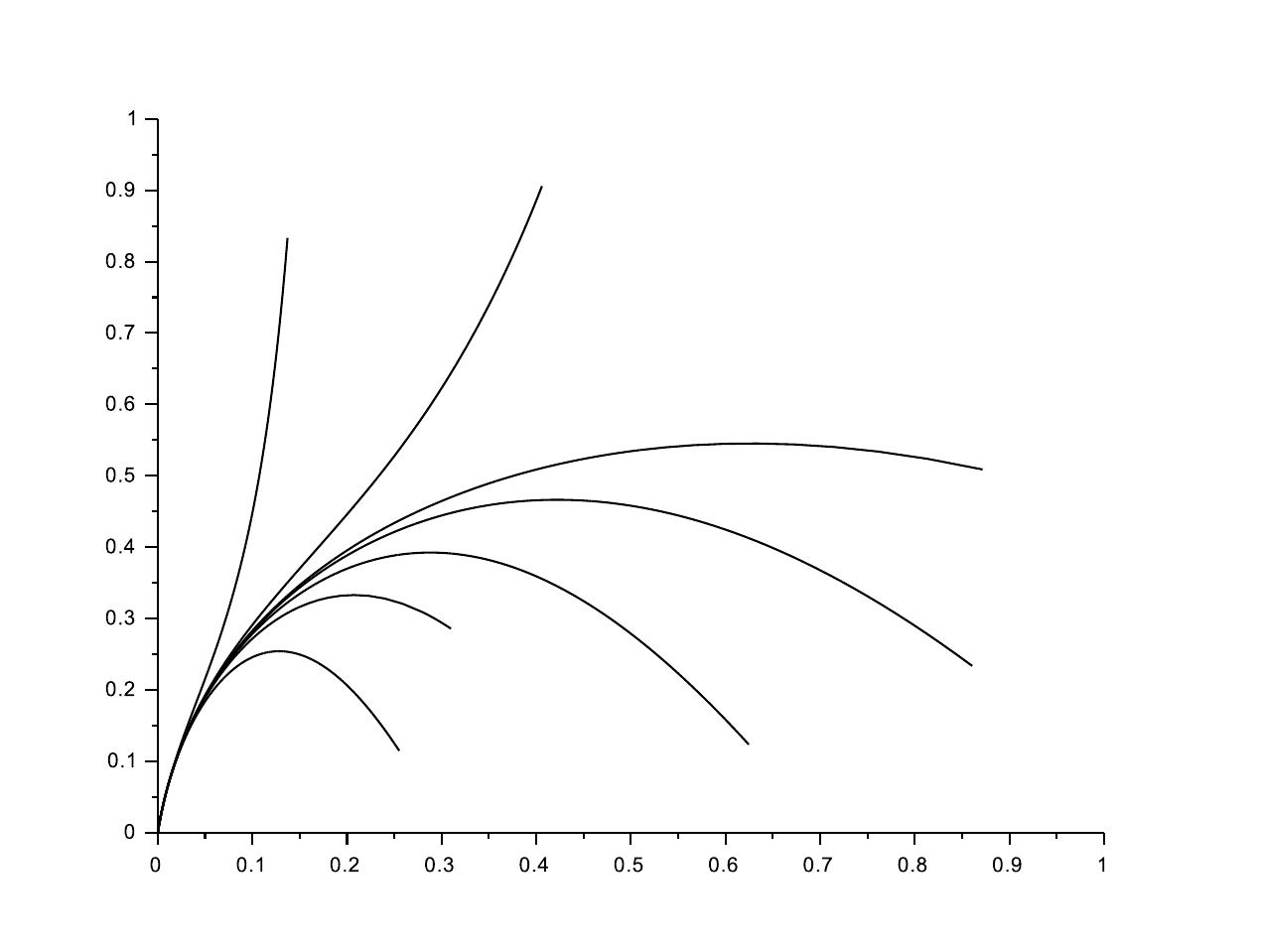}
   \end{minipage}
\caption{Example \ref{AstaCoins}, Phase portrait of $F^0$ and $F^1$}

\end{figure}
\begin{figure}[!h]
\centering
\includegraphics[scale=0.4]{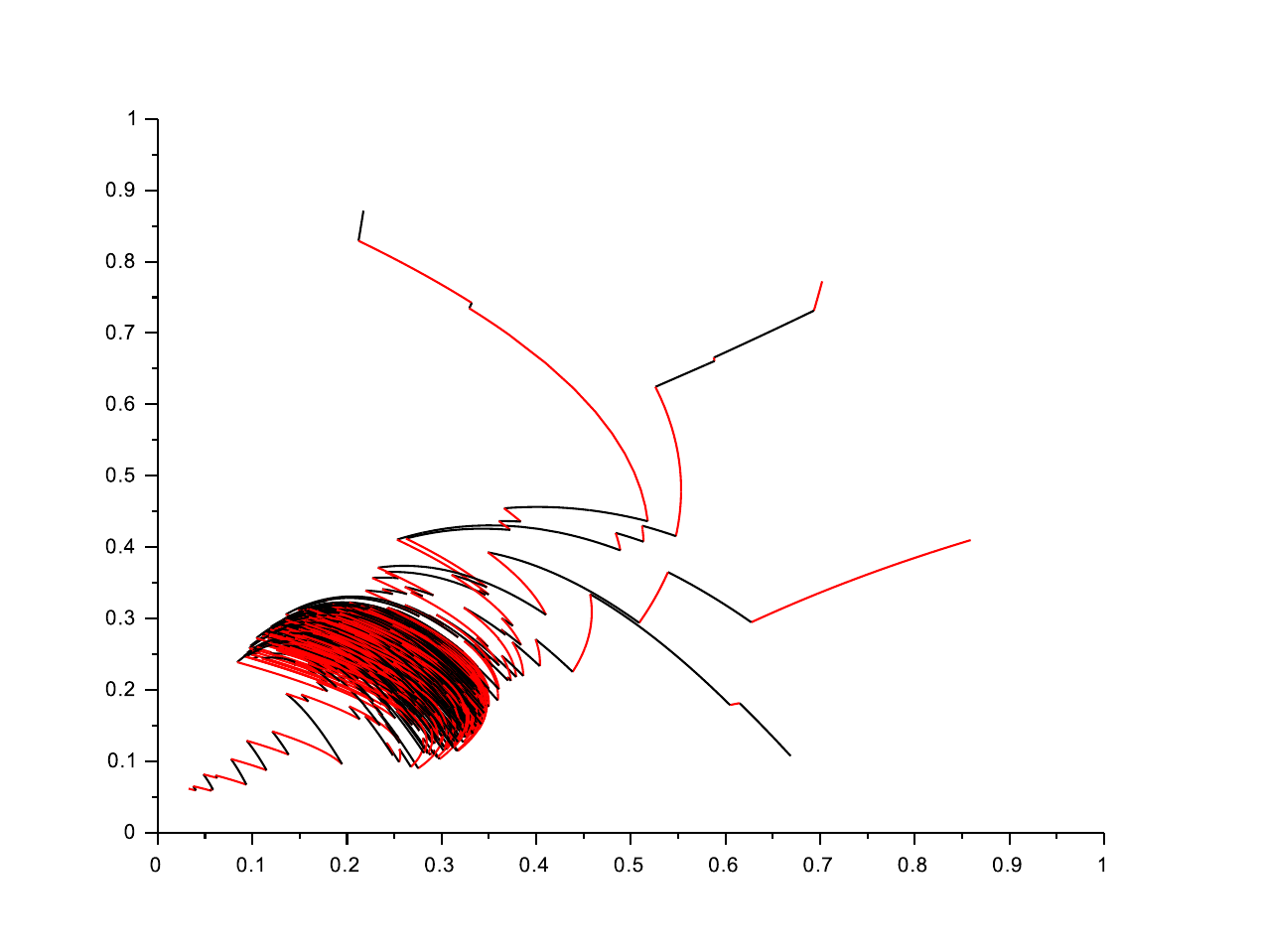}
\caption{\label{fig:switch} Example \ref{AstaCoins}, some trajectories of $(X_t)$ for $\beta = 20$}
\end{figure}
  \label{AstaCoins}
  }
  \eex

\brem {\rm 
In \cite{BHS18}, we show that the previous example can be generalised in the following way. Assume that $F^0$ and $F^1$ are two epidemic vector fields in dimension $2$ such that 
\begin{enumerate}
\item  $\lambda(A^0)<0$ and $\lambda(A^1)<0$,
\item There exists $s \in (0,1)$ such that $\lambda(A^s)>0$, where $A^s = s A^1 + (1-s) A^0$.
\end{enumerate}
Then, \cite[Lemma 3.7]{BHS18} show that there exists an accessible point at which the weak bracket condition holds. Moreover, since $\lambda(A^s)>0$, Theorem \ref{th:LajYorkeGen} implies that condition \textbf{(ii)} of Theorem \ref{th:persist3} is satisfied. Thus, by this theorem, we can conclude that there is convergence in total variation to a unique invariant probability measure provided $\lambda_1 > 0$. This happens for example with switching rates of the form $$a_{0,1} = s \beta,  a_{1,0} = (1-s) \beta, a_{0,0}=a_{1,1}=0.$$ for $\beta$ large enough (by Corollary \ref{cor:moyenne}.)}
\erem

  \brem
  \label{rem:FMG}{\rm
 In the preceding example, the matrices $A^i$ are Metzler and Hurwitz but   $\lambda_1 > 0$   because  the convex hull of the $\{A^i\}$ contains a non Hurwitz matrix. This leads to the natural question of finding examples for which:
 \begin{quote}  $\lambda_1 > 0$ and every matrix  in the convex hull of the $\{A^i\}$ is Hurwitz. \end{quote}
For arbitrary (i.e non Metzler)  matrices, such and example has been given in dimension 2 in \cite{Lawley&matt&Reed} and more recently in ~\cite{lagasquie16}.

Now, if we restrain ourselves to Metzler matrices, a result from Gurvits, Shorten and Mason (\cite[Theorem 3.2]{gurvits07})  proves that, in dimension 2, when every matrix in the convex hull is Hurwitz, then  $0$ is globally asymptotically stable for any deterministic switching between the linear systems. In particular, this  implies that $\lambda_1$ cannot be positive.

 However, they show that it is possible  in some  higher dimension  to construct  an example where all the matrices in the convex hull are Hurwitz, and for which there exists a periodic switching such that the linear system explodes. Later, an explicit example in dimension 3 was given by Fainshil, Margaliot and Chiganski \cite{onpls}. Precisely, consider the matrices

 $$A^0 = \begin{pmatrix}
   -1 & 0 & 0\\
   10 & -1 & 0\\
   0 & 0 & -10
\end{pmatrix}, \; A^1 = \begin{pmatrix}
   -10 & 0 & 10\\
   0 & -10 & 0\\
   0 & 10 & -1
\end{pmatrix}.$$
It is shown in \cite{onpls} that every convex combination of  $A^0$ and $A^1$ is Hurwitz, and yet a switch of period 1 between  $A^0$ and $A^1$ yields an explosion. Some simulations made on Scilab (see Figure \ref{fig:explo}) let us think that this result is still true for a random switching, with rates $$a_{0,1} = a_{1,0} = \beta, a_{0,0}=a_{1,1}=0.$$ Here $\beta$ has to be chosen neither too small nor too big.
 Using the formula $$\lim_{t \to \infty} \mathbb{E}(\frac{1}{t} \int_0^t \langle A^{J_s} \Theta_s, \Theta_s \rangle ds ) = \lambda_1(\beta),$$ and Monte-Carlo simulations we  can estimate numerically  $\lambda_1(\beta).$  The results  are plotted in Figure \ref{fig:lambda}  and show (although we didn't prove it) that $\lambda_1 > 0$ for $3 \leq \beta \leq 30,$  providing a positive answer to the question raised at the beginning of the remark.

\begin{figure}[!h]
\centering
\includegraphics[scale=0.2]{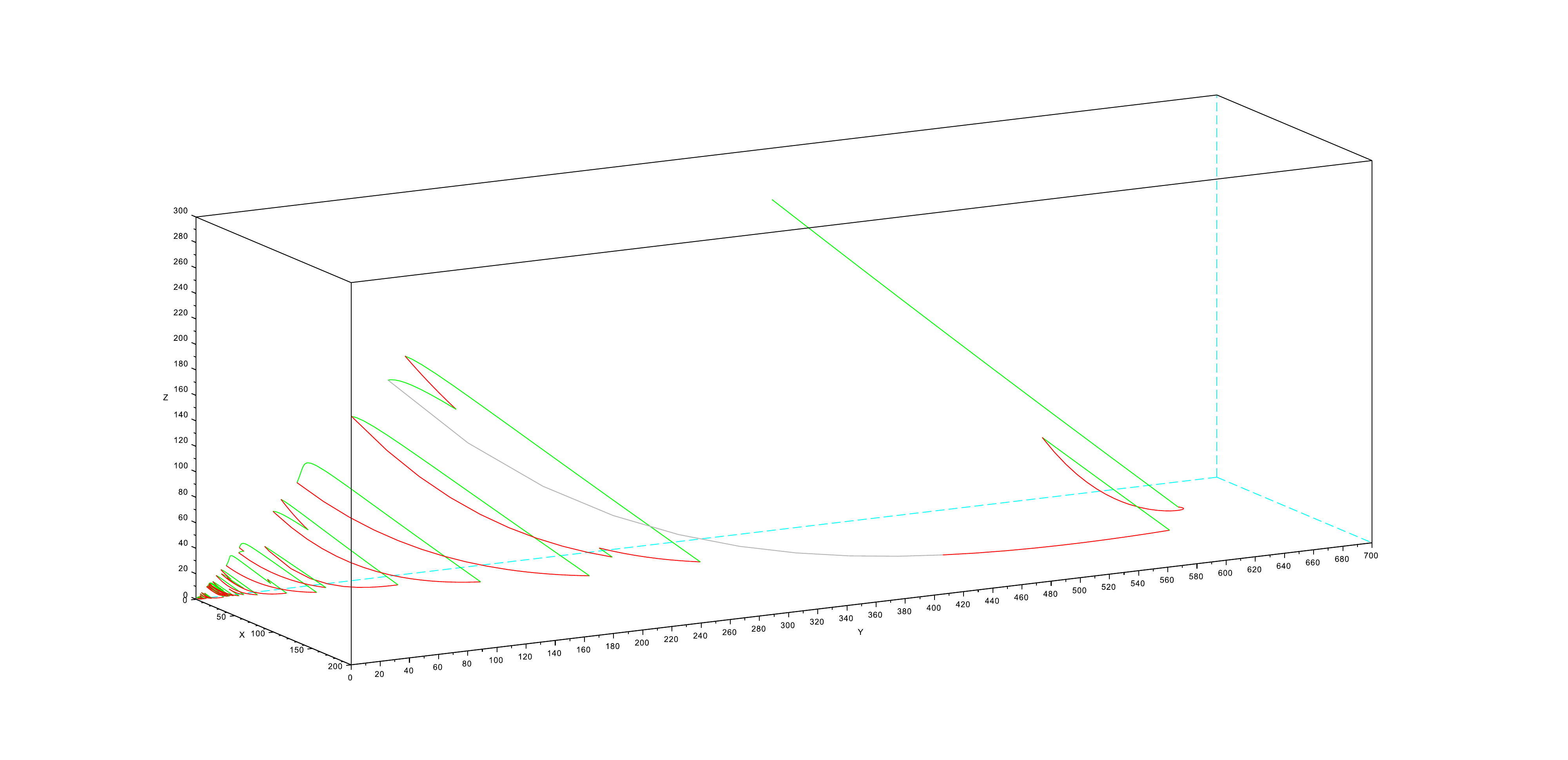}
\caption{\label{fig:explo} Simulation of $Y_t$ for $\beta = 10$.}
\end{figure}

\begin{figure}[!h]
\centering
\includegraphics[scale=0.2]{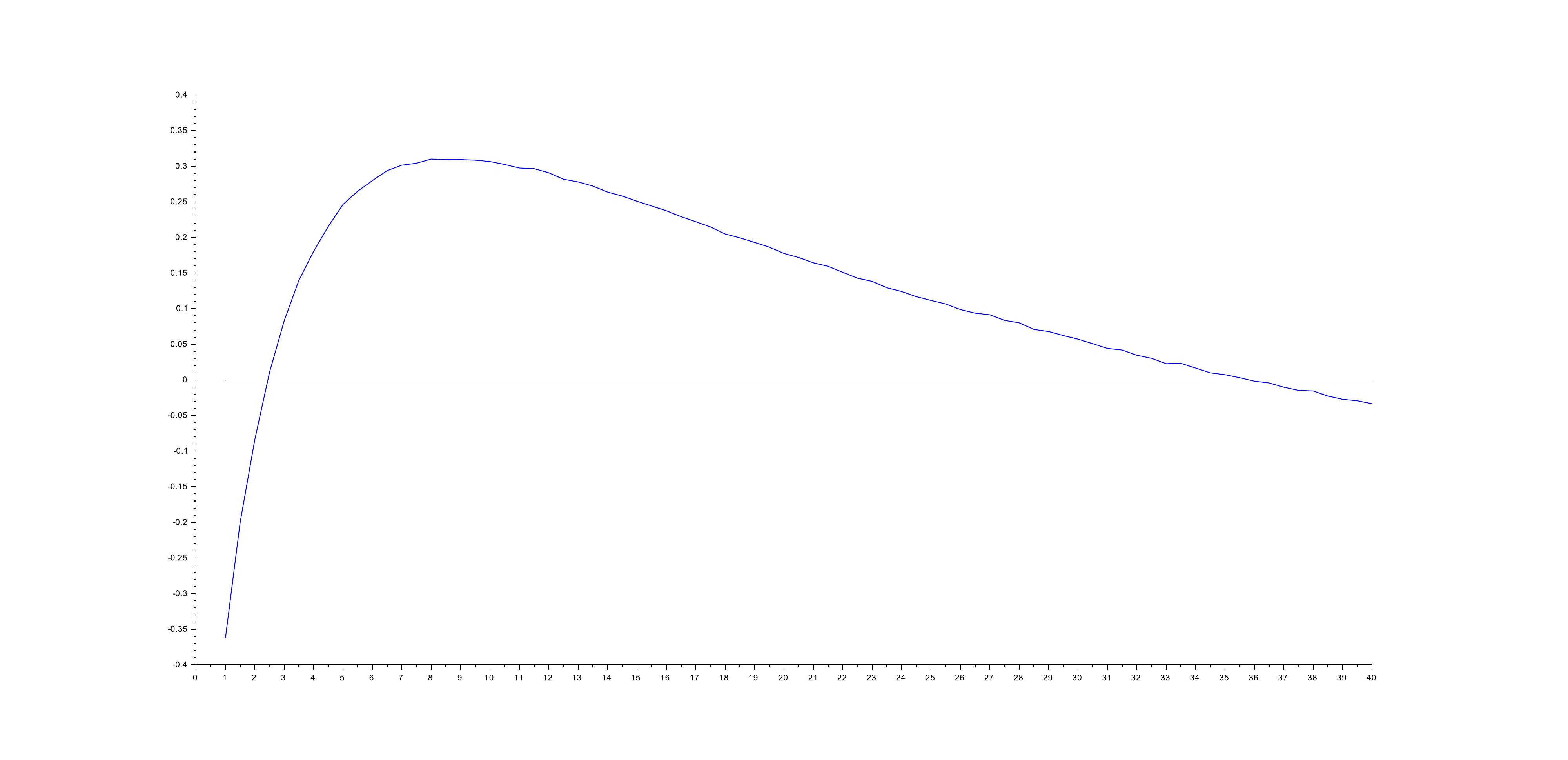}
\caption{\label{fig:lambda} Approximation of $\lambda_1(\beta)$ by Monte-Carlo method }
\end{figure}
}
\erem

\bex[Fluctuations may promote infection, continued]\label{ex:LY3d}{\rm Remark \ref{rem:FMG} can be used to produce two  Lajmanovich-Yorke vector fields $F^0, F^1$ on $[0,1]^3$
 such that
\bdes
\iti For all $0 \leq t \leq 1,$ the disease free equilibrium  is a global attractor of  the vector field  $F^t = (1-t) F^0 + t F^1;$
\itii A random switching between the dynamics leads to the persistence of the disease.
\edes
Observe that $F^t$ is the  Lajmanovich-Yorke vector field with infection matrix $C^t = (1-t) C^0 + t C^1$ and cure rate vector $D^t = (1-t) D^0 + t D^1$

To do so, one  just has to choose  $C^0,C^1,D^0,D^1$ in such way that $A^i=C^i-D^i$. For the simulation given here, we have chosen
$$ D^0  = \begin{pmatrix}
   11  \\
    11 \\
    20
\end{pmatrix},$$ and $$ D^1 = \begin{pmatrix}
   20  \\
   20 \\
   11
\end{pmatrix}.$$

When (see Figure \ref{fig:lambda}) $\beta$ is such that  $\lambda_1 > 0$, then by Theorem \ref{th:cvexpo} below, $Z$ admits a unique invariant measure $\Pi$ on $M^* \times E$. Moreover by Theorem \ref{th:persist1}, there exists $\theta > 0$ such that
$$\sum_{i \in E} \int \|x\|^{-\theta} \Pi^i(dx) < \infty.$$
Figure \ref{fig:LY3d} and \ref{fig:norm} illustrate this persistence of the infection. In figure \ref{fig:norm} , we have plotted $\|X_t\|_1 = X^1_t + X^2_t+ X^3_t$.

\begin{figure}[!h]
\centering
\includegraphics[scale=0.4]{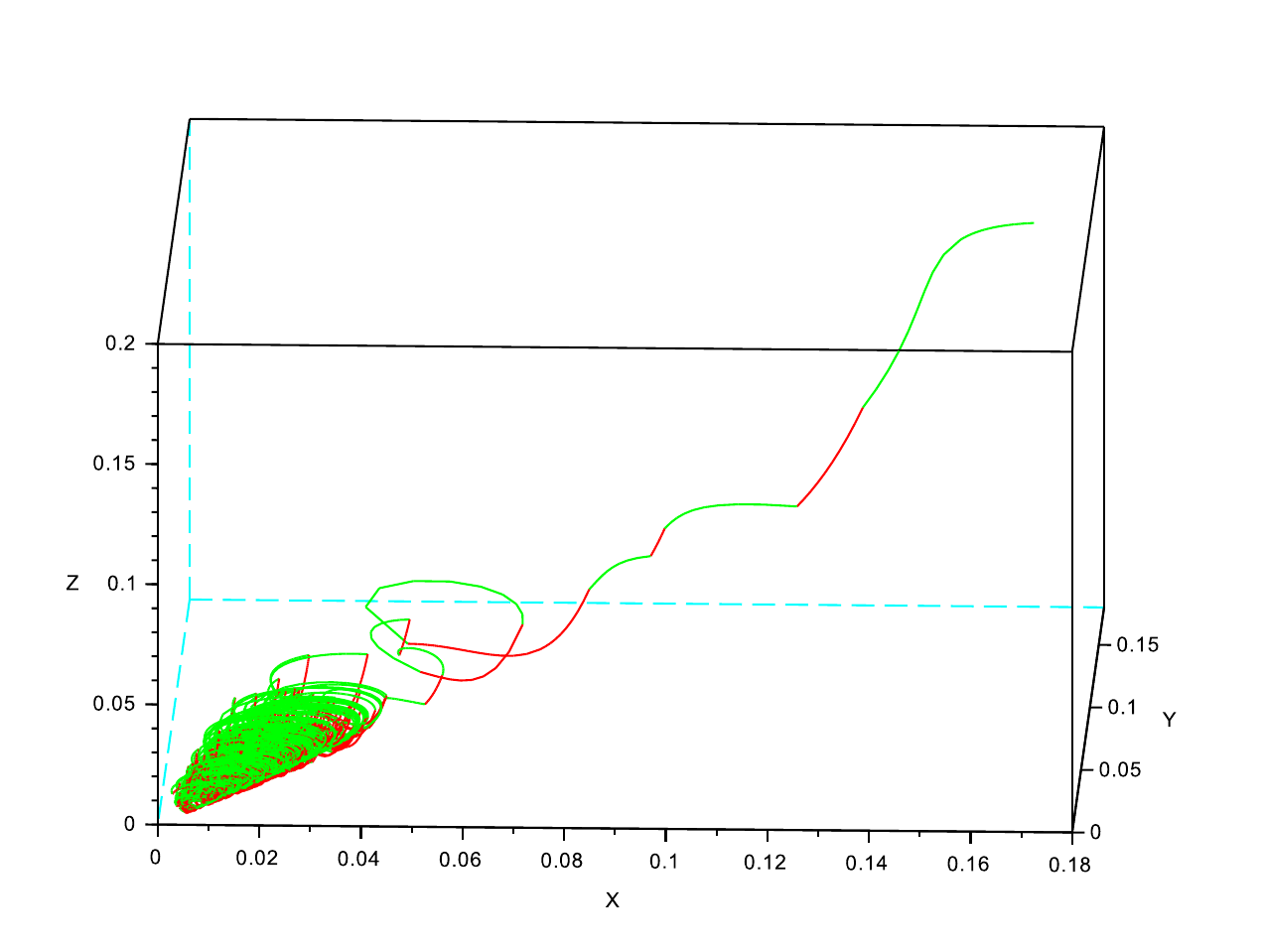}
\caption{\label{fig:LY3d} Example \ref{ex:LY3d} : Simulation of $X_t$ for $\beta = 10$.}
\end{figure}

\begin{figure}[!h]
\centering
\includegraphics[scale=0.4]{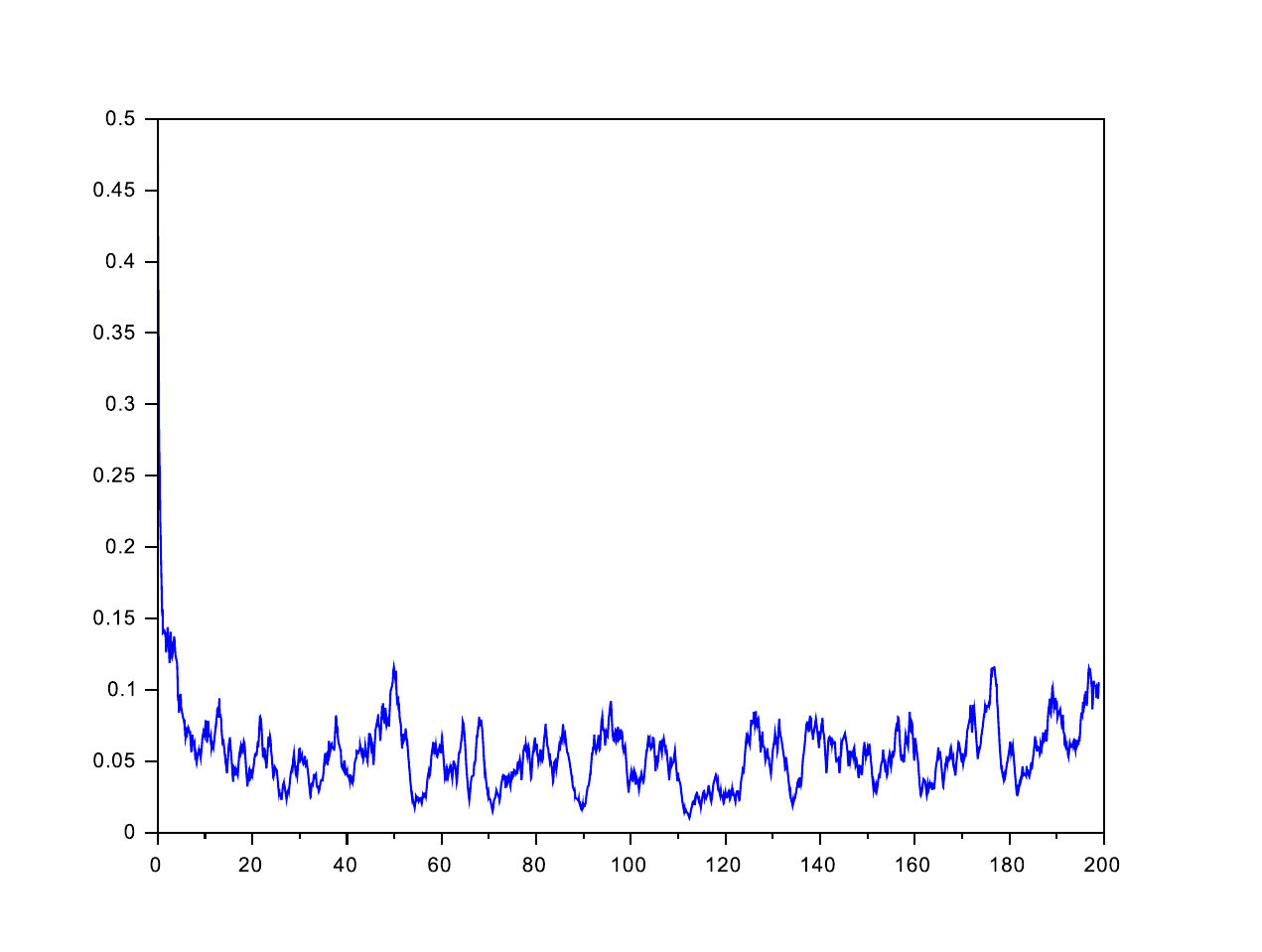}
\caption{\label{fig:norm} Example \ref{ex:LY3d} : Simulation of $\|X_t\|_1$ for $\beta = 10$. }
\end{figure}

}\eex

\subsection{Exponential convergence without bracket condition}
Throughout this section, we assume that the vector fields $F^i$ are epidemic and that the jump rates are constant. Recall (see proof of Theorem~\ref{th:EpidemicExtinct}) that this implies that for all $\omega \in \Omega$ and $t > 0$, $\Psi(t,\omega)$ is monotone and strongly subhomegeneous. A very useful consequence of this fact is the strict nonexpansivity of $\Psi(t,\omega)$ on $\RR^d_{++}$ with respect to the Birkhoff part metric $p$, the definition of which is recalled below. Now if we assume that $\lambda_1 > 0$, we have a Lyapunov function and nonexpansivity, so we might expect uniqueness of the invariant measure on $[0,1]^d \setminus \{0 \} \times E$ and convergence in law of $(Z_t)$ towards it. Here we prove that this is indeed the case, and even that we have an exponential rate of convergence towards this invariant measure for a certain Wasserstein distance, thanks to a weak form of Harris' theorem given by Hairer, Mattingly and Scheutzow~\cite{hms}. But before to do so, we explain briefly why we cannot expect to have convergence in total variation  without additional assumptions with the following simple example :
\bex{\rm
\label{ex:nobra}
Suppose $d = 2, E = \{0,1\}.$ Let  $F^0, F^1$ be the {\em Lajmanovich-Yorke vector fields} respectively given by  $$C^0 = \begin{pmatrix}
   1 & 3 \\
   2 & 4
\end{pmatrix}, \: D^0  = \begin{pmatrix}
   2 \\
    3
\end{pmatrix},$$ and $$C^1 = \begin{pmatrix}
   6 & 2 \\
   7 & 3
\end{pmatrix}, \:   D^1 = \begin{pmatrix}
   4 \\
    5
\end{pmatrix}.$$
One can easily check that the point $x^*=(1/2,1/2)$ is a common equilibrium of $F^1$ and $F^2$. In particular, $\Pi=\delta_{x^*} \otimes (\delta_0+ \delta_1)/2$ is an invariant probability of $Z.$  Moreover, for all $x \neq x^*$, $i \in E$ and $t \geq 0$, one has  $\mathbb{P}^Z_{x,i}( Z_t \in \{x^*\} \times E) = 0$ so $\| \delta_{x,i}P_t^Z - \Pi\|_{TV} = 1$ for all $t \geq 0$. Now let us quickly show that $X_t$ converges almost surely exponentially fast to $x^*$, for all switching rates. Let $\lambda_1(0) = \lambda_1$ (respectively $\lambda_1(x^*)$) denote   the top Lyapunov exponent of the linearized system at the origin (respectively at $x^*$). By Proposition \ref{th:Perron} this exponent coincides with the unique average growth rate of the corresponding linearized system. We claim that $\lambda_1(0) > 0$ and $\lambda_1(x^*) < 0.$ The first inequality follows from  the Kolotilina-type lower estimate for the  top Lyapunov exponent mentioned in Remark \ref{rem:Mierc} due to Mierczy{\'n}ski (\cite[Theorem 1.3]{Mier13}). In our setting, this estimate ensures that
$$\lambda_1(0) \geq \frac{1}{2} \sum_i p_i Tr(A^i) +  \sum_i p_i \sqrt{A_{12}^i A_{21}^i},$$
which is positive because $Tr(A^0)=Tr(A^1)=0$ and the other terms are positive.  Let $B^i = DF^i(x^*).$ Then the second estimate follows from  Lemma \ref{lem:roughestimate} because one can easily check that  $\lambda_{\max} (B^1+(B^1)^T) \leq\lambda_{\max} (B^0+(B^0)^T) < 0$. So applying Theorem \ref{th:extinct}, we have a neighborhood  ${\cal U}$ of $x^*$  and $\eta > 0$ such that for all $x \in {\cal U}$ and $i \in E$
\beq
\label{inqnobra1}
\mathbb{P}^Z_{x,i}( \limsup_{t \rar \infty} \frac{1}{t} \log(\|X_t-x^*\|) \leq \frac{\lambda_1(x^*)}{2}) \geq \eta.
\eeq
On the other hand, because $\lambda_1(0) > 0$, there exists by Theorem \ref{th:persist1} $\varepsilon > 0$ such that for all $x \neq 0$,
\beq
\label{inqnobra2}
 \mathbb{P}^Z_{x,i}( \tau < \infty) = 1,
\eeq
where $\tau = \inf \{ t \geq 0 \: : \: \| X_t |\ \geq \varepsilon \}$. Finally, because $x^*$ is a linear stable equilibrium for $F^0$ with basin of attraction contains $M^*$, one can show that there exists a constant $c > 0$ such that for all $x \in M$ with $\|x\| \geq \varepsilon$,
\beq
\label{inqnobra3}
  \mathbb{P}^Z_{x,i} ( Z_t \in \mathcal{U} \times E) \geq c.
\eeq
Combining \eqref{inqnobra1},  \eqref{inqnobra2}, \eqref{inqnobra3} and the Markov property implies that
$$\mathbb{P}^Z_{x,i}( \limsup_{t \rar \infty} \frac{1}{t} \log(\|X_t-x^*\|) \leq \lambda_1(x^*))=1,$$
for all $(x,i) \in M^* \times E$ (see \cite[Theorem 3.1]{BL16} for details on a very similar proof).
}
\eex

Before stating our theorem, recall the definition of the Wasserstein distance. Let $\mathcal{Y}$ be a Polish space, and $d$ be a distance-like function on  $\mathcal{Y}$. That is $d$ satisfies the axioms of a distance, except for the triangle inequality. Then the \textit{Wasserstein distance} associated to $d$ is defined for every $\mu, \nu \in {\cal P}(\mathcal{Y})$ by

$$ {\cal W}_d(\mu,\nu) = \inf_{\pi \in C(\mu,\nu)} \int_{\mathcal{X}^2} d(x,y) \mathrm{d}\pi(x,y),$$

where $C(\mu,\nu)$ is the set of all the coupling of $\mu$ and $\nu$. When $d$ is a distance, so is ${\cal W}_d$, and in every case, ${\cal W}_d(\mu,\nu) =0$ if and only if $\mu = \nu$.

Set $\mathcal{Y}= [0,1]^d \setminus \{0\} \times E$.




\bthm
Assume the $F^i$ are epidemic vector fields, $(a_{ij})$ are constant and $\lambda_1 >0$. Then there exists a distance-like function $\tilde{d}$, $t_0 \geq 0$ and $r >0$, such that,

\bdes
\iti for all $t \geq t_0$, for all $\mu, \nu \in {\cal P}(\mathcal{Y})$,
$$ {\cal W}_{\tilde{d}}( \mu P_t^Z, \nu P_t^Z) \leq e^{- r t} {\cal W}_{\tilde{d}}( \mu, \nu ).$$
\itii $(P_t^Z)$ has a  unique invariant measure $\Pi$   on $\mathcal{Y}$, and for all $\mu \in {\cal P}(\mathcal{Y})$,
$$ {\cal W}_{\tilde{d}}( \mu P_t^Z, \Pi) \leq e^{- r t} {\cal W}_{\tilde{d}}( \mu, \Pi ).$$
\edes

\label{th:cvexpo}
\ethm

\section{Proofs of Theorems \ref{th:extinct}--\ref{th:persist3} : A {\em stochastic persistence} approach}
\label{sec:mainproofs}
As indicated in the introduction, the proofs will be deduced from the qualitative properties of PDMPs  combined with general results on {\em stochastic persistence} proved in \cite{Ben14} along the lines of the seminal results obtained by Schreiber, Hofbauer and their co-authors for deterministic systems.
\subsection{An abstract stochastic persistence result}
The results in \cite{Ben14} concern certain Markov processes on a (possibly) non compact metric space satisfying a weak version of the Feller property. Here for simplicity we shall state a simpler version of these results tailored for Feller processes on a compact space.

Let $\XX$ be a  {\bf compact} metric  space and $\tilde{Z}$ a cad-lag Markov process on $\XX.$ To shorten notation we write
 $\mathbb{P}_x, \mathbb{P}_{\mu}, (P_t)_{t \geq 0}, {\cal P}_{inv}, {\cal P}_{erg}$  in place of  $\mathbb{P}^{\tilde{Z}}_x, \mathbb{P}^{\tilde{Z}}_{\mu}, (P^{\tilde{Z}}_t)_{t \geq 0}, {\cal P}_{inv}^{\tilde{Z}}, {\cal P}_{erg}^{\tilde{Z}}$. We let
  $$\Pi_t = \frac{1}{t} \int_0^t \delta_{\tilde{Z}_s} ds$$ denote the empirical occupation measure of $\tilde{Z}.$ We let $C(\XX)$ denotes the space of real valued continuous functions on $\XX$ equipped with the uniform norm $\|f\| = \sup_{x \in \XX} |f(x)|.$

We assume that $(P_t)_{t \geq 0}$ is {\em Feller}. That is
\bdes
\ita For all $t \geq 0$ $P_t$ maps $C(\XX)$ into itself,
\itb  For all $f \in C(\XX) \lim_{t \rar 0} \|P_t f - f\| = 0.$
\edes
We let $\L$ denote the {\em infinitesimal generator} of $(P_t)$ and $\D$ its domain. Recall that $\D$ is defined as the set of $f \in C(\XX)$
such that $\frac{1}{t} (P_t f - f)$ converges in $C(\XX),$ and, for such an $f,$ $\L f$ denotes the limit.
We let $\D^2 \subset \D$ denote the set of $f \in \D$ such that $f^2 \in \D.$ For $f \in \D^2$ the {\em Carr\'e du champ} of $f$ is defined as
\beq
\label{defGamm}
\Gamma(f) = \L f^2 - 2 f \L f.
\eeq
 We assume that
\begin{hypothesis}
\label{hyp:extinct}
 there exists a non empty compact set $\XX_0 \subset \XX$  called the {\em extinction set} which is invariant under $(P_t)_{t \geq 0}.$ That is
$$P_t \1_{\XX_0} = \1_{\XX_0}$$ where $\1_{\XX_0}$ stands for the indicator function of $\XX_0.$
\end{hypothesis} We set $$\XX_+ = \XX \setminus \XX_0,$$
 ${\cal P}_{inv}(\XX_+) = {\cal P}_{inv} \cap {\cal P}(\XX_+), {\cal P}_{inv}(\XX_0) = {\cal P}_{inv} \cap {\cal P}(\XX_0)$ etc.

{\em Extinction}  of $\tilde{Z}$ amounts to say that trajectories of $(\tilde{Z}_t)$ converge almost surely to $\XX_0.$ Let $\XX_0^{\eps}$ be the $\eps$-neighborhood of $\XX_0.$ Using a terminology borrowed to  Schreiber \cite{Sch12} and Chesson \cite{C82}, we say that $\tilde{Z}$ is  {\em stochastically persistent} (or {\em almost surely  persistent}), respectively {\em persistent in probability},
provided $$\lim_{\eps \rar 0 } \limsup_{t \rar \infty} \Pi_t(\XX_0^{\eps})  = 0$$ $\mathbb{P}_x$ almost surely for all $x \in \XX_+.$ Respectively
 $$\lim_{\eps \rar 0 } \limsup_{t \rar \infty} \mathbb{P}_x( Z_t \in \XX_0^{\eps}) = 0$$ for all $x \in \XX_+.$

General criteria ensuring extinction or persistence are given by the existence of a suitable {\em average Lyapounov function} $V$ as defined now.

In addition to hypothesis \ref{hyp:extinct}  we  assume that
\begin{hypothesis}
\label{hyp:lyap}
There exist  continuous maps $V : \XX_+ \mapsto \RR^+$ and $H : \XX \mapsto \RR$ enjoying the following properties :
\bdes
\ita For all compact $K \subset \XX_+$ there exists $V_K \in \D^2$ with  $V|_K = V_K|_K$ and  $(\L V_K)|_K = H|_K;$
\itb $\sup_{\{K : K \subset \XX_+, K \mbox{ compact } \}} \|\Gamma(V_K)|_K\|   < \infty;$
\itc $\lim_{x \rar \XX_0} V(x) = \infty;$
\itd Jumps of $V(\tilde{Z}_t)$ are bounded : $\exists \Delta > 0$ such that $|V(\tilde{Z}_t) - V(\tilde{Z}_{t-})| \leq \Delta;$

\edes
\end{hypothesis}
 Let ${\cal P}_{erg}(\XX_0) = {\cal P}_{erg} \cap {\cal P}(\XX_0).$ Define the $H$-exponents of the processes as
 $$ \Lambda^+(H)  = - \inf_{\mu \in {\cal P}_{erg}(\XX_0)} \mu H  \mbox{ and } \Lambda^-(H) = - \sup_{\mu \in {\cal P}_{erg}(\XX_0)} \mu H.$$ We call the process {\em $H$-persistent} if $\Lambda^-(H) > 0$ and
{\em $H$-nonpersistent} if $\Lambda^+(H) < 0.$

By the Ergodic decomposition theorem, note that $\Lambda^-(H) > 0$ (respectively $\Lambda^+(H) > 0$) if and only if $\mu H < 0$ (respectively $ > 0$) for all $\mu \in {\cal P}_{erg}(\XX_0).$
\\

We say that $A \subset \XX$ is  {\em accessible} from  $B \subset \XX$ if for every neighborhood $U$ of $A$ and $x \in B$ there exists $t \geq 0$ such that $P_t \1_U (x) > 0.$

We call a point $p \in \XX$ a {\em Doeblin point} provided there exists a neighborhood $U$ of $p$, a bounded (positive) measure $\nu$ on $\XX$ and some  number $s > 0$ such that
$$\delta_x P_s \geq \nu$$ for all $x \in U.$ The following theorem is a consequence of Theorems 4.4 and 4.10 and Proposition 8.2 in \cite{Ben14}.
\bthm
\label{th:Hpers}Suppose that the process is $H$-persistent. Then
\bdes
\iti The process is stochastically persistent. In particular, for all $x \in \XX_+,$ $\mathbb{P}_x$ almost surely, every limit point of $\{\Pi_t\}$ lies in
${\cal P}_{inv}(\XX_+) = {\cal P}_{inv} \cap {\cal P}(\XX_+).$
\itii There exist $0 < \rho < 1$ and positive constants $\theta > 0, K > 0, T$ such that
$$P_T (e^{\theta V}) \leq \rho e^{\theta V} + K;$$
\itiii  Let $\eps > 0$ and $\tau^{\eps}$ be the stopping time defined by
$$\tau^{\eps} = \inf \{t \geq 0 : \: \tilde{Z}_t\in \XX_0^{ \eps} \}.$$ Then there exists $\eps >0$ such that  for all  $1 < b < \frac{1}{\rho}$, there exists $c > 0$  such that for all $x \in \XX^+$
$$\mathbb{E}_x(b^{\tau}) \leq c (1 + e^{\theta V(x)});$$
\itiv If, furthermore, there exists a Doeblin point $x \in \XX_+$ accessible from $\XX_+$ then ${\cal P}_{inv}(\XX_+)$ reduces to a single measure $\Pi$ and
for all $x \in \XX_+$
$$\|\delta_x P_t - \Pi \|_{TV} \leq const.(1 + e^{\theta V(x)}) e^{-\kappa t}$$
for some $\kappa > 0.$
\edes
\ethm

The next result is a general extinction result.
\bthm
\label{th:Hnonpers}
Suppose that the process is $H$-nonpersistent.  Then
\bdes
\iti For all $0 < \alpha < - \Lambda^+(H)$, there exists a neighborhood $U$ of $\XX_0$ and $\eta >0$ such that
$$\mathbb{P}_x(\liminf_{t \rar \infty} \frac{V(\tilde{Z}_t)}{t} \geq \alpha)\geq \eta$$ for all $x \in U;$
\itii If furthermore $\XX_0$ is accessible from $\XX$
$$\mathbb{P}_x(\liminf_{t \rar \infty} \frac{V(\tilde{Z}_t)}{t} \geq - \Lambda^+(H)) = 1$$ for all $x \in \XX.$
\edes
\ethm

\prf Since the proof is very similar to the one given in \cite[Theorem 3.1]{BL16}, we only give a sketch of it. Let $0 < \alpha < - \Lambda^+(H)$. The proofs of Propositions 8.2 and 8.3 in \cite{Ben14} (see also \cite[Lemma 3.5]{BL16})  adapt verbatim in the nonpersistent case to prove that there exist $T > 0$, $\theta > 0$, $\eps > 0$ and $0 < \rho < 1$ such that, for all $z \in \XX_0^{\eps} \setminus \XX_0$, 

\begin{enumerate}
\item[(i)] $P_T V (z) - V(z) \geq \alpha T$,
\item[(ii)] $P_T e^{- \theta V}(z) \leq \rho e^{-\theta V(z)}$.
\end{enumerate}
Here and throughout this proof, $\XX_0^{\eps} = \{ z \in \XX_+ \: : \: V(z) > - \log( \eps )\} \cup \XX_0$. We set $\tau_{\eps} = \inf \{ k \geq 0 \: : \: \tilde{Z}_{kT} \notin \XX_0^{\eps} \}$. We claim that 

\begin{enumerate}
\item There exists $\eta > 0$ such that for all $z \in \XX_0^{\eps/2}$, $\mathbb{P}_z( \tau_{\eps} = \infty) \geq \eta$;
\item On the event $\{ \tau_{\eps} = \infty \}$, and for all $z \in \XX_0^{\eps/2}$,  $\liminf_{t \rar \infty} \frac{V(\tilde{Z}_t)}{t} \geq \alpha$.
\end{enumerate}
In particular, this implies point \textbf{(i)} of the Theorem with $U = \XX_0^{\eps/2}.$ Point \textbf{(ii)} easily follows by Markov property. We prove the first claim. We set for $k \geq 0$, $W_k  = e^{-\theta V(\tilde{Z}_{kT})}$. Due to point (ii) above, $(W_{k \wedge \tau_{\eps}})_{k \geq 0}$ is a supermartingale. In particular, for all $z \in \XX_0^{\eps/2} \setminus \XX_0$, 

$$\mathbb{E}_z ( W_{k \wedge \tau_{\eps}} \1_{\tau_{\eps} < \infty} ) \leq e^{-V(z)} \leq \left(\frac{\eps}{2}\right)^{\theta}.$$ By dominated convergence, this gives $$\eps^{\theta} \mathbb{P}_z( \tau_{\eps} < \infty ) \leq \left(\frac{\eps}{2}\right)^{\theta},$$ which proves the first point with $\eta =  1-2^{-\theta}$. We now prove the second claim. We set $M_n = \sum_{k=1}^n P_T V( \tilde{Z}_{(k-1)T}) - V(\tilde{Z}_{kT})$. The sequence $(M_n)_{n \geq 1}$ is a martingale, and on the event $\{ \tau_{\eps} = \infty \}$ and for all $z \in \XX_0^{\eps/2} \setminus \XX_0$, 

$$\frac{M_n}{n} \geq \alpha - \frac{V(\tilde{Z}_{nT})}{n}.$$ Hence the strong law of large numbers for martingales implies that, on the event $\{ \tau_{\eps} = \infty \}$ and for all $z \in \XX_0^{\eps/2} \setminus \XX_0$,  
\beq 
\label{eqn:n}
\liminf_{n \to \infty} \frac{V(\tilde{Z}_{nT})}{T} \geq \alpha.
\eeq Now, Lemma 7.4 in \cite{Ben14} implies that for all $z \in \XX_+$, the process $$ M_t^V = V(\tilde{Z}_t) - V(z) - \int_0^t H(\tilde{Z}_s) d s$$ is a martingale such that, almost surely, $\lim_{t \to \infty} \frac{M_t^V}{t} =0$. Since $H$ is bounded, this implies that for all $t \in [0,T]$, $$\lim_{n \to \infty}\frac{V(\tilde{Z}_{nT + t}) - V(\tilde{Z}_t)}{n}=0.$$ This, together with \eqref{eqn:n} proves the second claim. \qed

\subsection{Proofs of Theorems \ref{th:extinct}--\ref{th:noncompactglobal}}
 In order to apply  the results of the previous section we rewrite the dynamics of $Z = (X,I)$ in polar coordinates. Let $\Psi : M^*  \times E \to \RR_{+}^* \times S^{d-1} \times E$ be defined by $\Psi(x,i) = (\|x\|, \frac{x}{\|x\|}, i)$ and
  $$\mathcal{X}_+ = \Psi(M^* \times E).$$ Whenever $X_0 \in M^*,$ the process $\tilde{Z}_t = \Psi(Z_t)=(\rho_t, \Theta_t, I_t) \in \mathcal{X}_+$ satisfies the system
\beq
\label{dThetarhoNL}
 \left \{
\begin{array}{l}
 \frac{d \rho_t}{dt} =  \langle \Theta_t, \tilde{F}^{I_t}(\rho_t,\Theta_t) \rangle \rho_t \\ \\
 \frac{d\Theta_t}{dt} = \tilde{F}^{I_t}(\rho_t,\Theta_t)-\langle \Theta_t, \tilde{F}^{I_t}(\rho_t,\Theta_t) \rangle \Theta_t \\ \\
  \Pr(I_{t+s} = j | {\cal F}_t) = a_{i j}(\rho_t \Theta_t) s + o(s) \mbox{ for } i \neq j \mbox{ on } \{I_t = i \}
\end{array}
\right.
\eeq
where $$\tilde{F}^{i}(\rho,\theta)= \frac{F^i(\rho \theta)}{\rho}$$ for all $\rho >  0$ and $\theta \in S^{d-1}.$
By $C^2$ continuity of $F^i,$ the map $\tilde{F}^i$ extends to a $C^1$ map $\tilde{F}^i : \Rp \times S^{d-1} \mapsto \RR^d$ by setting
 $$\tilde{F}^{i}(0,\theta)= A^i \theta.$$
Thus, using this extension, (\ref{dThetarhoNL}) extends to the state space
 $$\mathcal{X} := \overline{\mathcal{X}_+} = \mathcal{X}_+ \cup \mathcal{X}_0$$
where $\mathcal{X}_0 = \{ 0 \} \times (S^{d-1} \cap C_M) \times E.$

This induces a PDMP  (still denoted $\tilde{Z}$) on $\mathcal{X},$ whose infinitesimal generator $\tilde{\mathcal{L}}$ acts on functions $f : \mathcal{X} \to \RR$ smooths in $(\rho, \theta)$ according to
\beq
\label{def:Lpolaire}
\tilde{\mathcal{L}} f(\rho, \theta, i) = \frac{\partial f^i}{\partial \rho}(\rho, \theta) \langle \theta, \tilde{F}^i(\rho,\theta) \rangle \rho + \langle \nabla_{\theta} f^i (\rho, \theta), \tilde{G}^i(\rho, \theta) \rangle + \sum_{j \in E} a_{ij}(\rho \theta) (f^j(\rho, \theta) - f^i(\rho, \theta)),
\eeq
where $\tilde{G}^i(\rho, \theta) = \tilde{F}^{i}(\rho,\theta)-\langle \theta, \tilde{F}^{i}(\rho,\theta) \rangle \theta$.
By  \cite[Proposition 2.1]{BMZIHP}, $\tilde{Z}$ is Feller. Moreover by equation \eqref{dThetarhoNL}, Hypothesis \ref{hyp:extinct} is verified. The following lemma gives $V$ and $H$ that fulfil Hypothesis \ref{hyp:lyap}.
\blem
For all $(\rho, \theta, i) \in \mathcal{X}$, set $H(\rho, \theta , i) = - \langle  \tilde{F}^i(\rho, \theta), \theta \rangle$, and for $\rho \neq 0$, $V(\rho, \theta , i)  =  -\log(\rho)$. Then $V$ and $H$ satisfy Hypothesis \ref{hyp:lyap}.
\elem
\prf
The definition of $\tilde{\mathcal{L}}$ and $V$ imply that $\tilde{\mathcal{L}} V (\rho, \theta, i) = H (\rho, \theta, i)$ for all $(\rho, \theta, i) \in \mathcal{X}_+.$  For all $K \subset \mathcal{X}_+$ compact, there exists $\eps > 0$  such that $\rho \geq \eps$ on $K.$ Let $\log_{\eps} : \RR \mapsto \RR$ be a smooth function coinciding with $\log$ on $[\eps, \infty[.$ Set $V_{K}(\rho, \theta,i) = - \log_{\eps}(\rho).$ Then \textbf{(a)} is satisfied, and because $V_K$ doesn't depend on $i$  $\Gamma(V_K) = 0$ so that \textbf{(b)} is also satisfied.   \textbf{(c)} and \textbf{(d)} are  clearly satisfied.
\qed

Now we link the $H$-exponents of $\tilde{Z}$ with the extremal average growth rates of $Z$ :
\blem
\label{lem:Hlambda}
With the notation of the previous sections,
$$\Lambda^+(H) =  \Lambda^+ \quad \text{and} \quad \Lambda^-(H) =  \Lambda^-.$$
In particular, $\tilde{Z}$ is $H$-persistent if and only if $\Lambda^- > 0$ and $H$-nonpersistent if and only if $ \Lambda^+ < 0$.
\elem
\prf
On $\mathcal{X}_0$, $\tilde{Z}_t = (0, \Theta_t, J_t)$
where $(\Theta_t, J_t)$  is the process given in Section \ref{sec:linear}.
 Now, $\langle A^i \theta , \theta \rangle = - H(0,\theta,i)$, and the result easily follows from the definitions of $\Lambda^{+/-}, \Lambda^{+/-}(H)$
\qed

Thanks to these lemmas and theorems of the previous sections, we can now prove our main results.

\prf \textbf{of Theorem~\ref{th:extinct}}
Here we assume $\Lambda^+ < 0$, thus by Lemma \ref{lem:Hlambda} $\tilde{Z}$ is $H$ - nonpersistent. Theorem \ref{th:Hnonpers} \textbf{(i)} then gives exactly the first part of Theorem~\ref{th:extinct} because $V(\tilde{Z_t}) = - \log( \rho_t) = - \log( \| X_t \|)$ for all $x \neq 0$.

Assume furthermore that $0$ is $F$ - accessible from $M$. By \cite[Proposition 3.14]{BMZIHP}, this implies that $\{0\} \times E$ is accessible from $M \times E$ for the process $Z$ and thus that $\mathcal{X}_0$ is accessible from $\mathcal{X}$ for the process $\tilde{Z}$. Then  Theorem \ref{th:Hnonpers} \textbf{(ii)} proves the second assertion of Theorem~\ref{th:extinct}.
\qed

To show the other theorems, we use the following lemma for which the proof is omitted. Here, $\varphi$ denotes $\Psi^{-1}$.
\blem
\label{lem:Pinvbij}
The map
\begin{displaymath}
\begin{array}{lrcl}
 \mathcal{P}_{inv}^{\tilde{Z}}(\mathcal{X}_+) & \longrightarrow &  \mathcal{P}_{inv}^Z( M^* \times E)  \\
     \Pi & \longmapsto & \Pi \circ \varphi^{-1}
\end{array}
\end{displaymath}
is a bijection. Moreover, for all $(x,i) \in M^* \times E$, and all $t \geq 0$
$$ \Pi_t^{x,i} = \tilde{\Pi}_t^{\Psi(x,i)} \circ \varphi^{-1}.$$
Thus, by bi-continuity of $\Psi$, $\Pi_t^{x,i}$ converges almost surely to some $\Pi$ if and only if $\tilde{\Pi}_t^{\Psi(x,i)}$ converges to $\Pi \circ \Psi^{-1}$.
\elem
\prf \textbf{of Theorem~\ref{th:persist1}}
Here we assume $\Lambda^- > 0$, thus by Lemma \ref{lem:Hlambda} $\tilde{Z}$ is $H$ - persistent. Then Theorem \ref{th:Hpers} \textbf{(i)} and Lemma \ref{lem:Pinvbij} imply \textbf{(i)} of Theorem \ref{th:persist1}. Moreover, by   Theorem \ref{th:Hpers} \textbf{(ii)}, we have for some positive $\theta, K, T$
$$\tilde{P}_T (e^{\theta V}) \leq \rho e^{\theta V} + K.$$
Let $\tilde{\mu} \in  \mathcal{P}_{inv}^{\tilde{Z}}(\mathcal{X}_+)$ and set $\tilde{W}= e^{\theta V}$. Then integrating the previous inequality against $\tilde{\mu}$ gives $\tilde{\mu} \tilde{W} \leq \rho \tilde{\mu} \tilde{W} + K$, thus
\beq
\label{inq:mu}
 \tilde{\mu} \tilde{W} \leq \frac{K}{1-\rho}.
\eeq
Now let $\mu \in \mathcal{P}_{inv}^Z( M^* \times E)  $ and set $W(x,i) = \| x \|^{-\theta}$. Then $\mu W = (\mu \circ \Psi^{-1} \circ \Psi ) W = (\mu \circ \Psi^{-1}) (W \circ  \Psi^{-1})$. By lemma \ref{lem:Pinvbij},  $\mu \circ \Psi^{-1} \in \mathcal{P}_{inv}^{\tilde{Z}}(\mathcal{X}_+)$, and because $W \circ \Psi^{-1} = \tilde{W}$, \eqref{inq:mu} proves \textbf{(ii)} of  Theorem \ref{th:persist1}. Point  \textbf{(iii)} is immediate from \textbf{(iii)} of Theorem \ref{th:Hpers}.
\qed

\prf \textbf{of Theorem~\ref{th:persist2}}
By Theorem \ref{th:persist1}, $\mathcal{P}_{inv}^Z( M^* \times E)$ is non-empty. So the weak bracket condition implies by \cite[Theorem 4.5]{BMZIHP} uniqueness of $\Pi$ and the absolute continuity. Moreover, for all $(x,i)\in M^* \times E$, $(\Pi_t^{x,i})_{t \geq 0}$ is tight and admits a unique limit point $\Pi$, so that $\Pi_t^{x,i}$ converges almost surely to $\Pi$.
\qed

\prf \textbf{of Theorem~\ref{th:persist3}}
Assume that  the weak bracket condition holds at a point $p$ that is $F$-accessible from $M^*$ and that condition \textbf{(i}) or \textbf{(ii)} of Theorem \ref{th:persist3} holds. Then \cite[Theorem 4.2]{BMZIHP} in case \textbf{(i)} (respectively \cite[Theorem 2.6]{BHS18} in case \textbf{(ii)}) implies that  for all $i$, $\Psi(p,i)$ (resp. $\Psi(e^{\star},i)$)  is a Doeblin point, which is accessible for the process $\tilde{Z}$ from $\mathcal{X}_+$. Thus by point \textbf{(iv)} of Theorem \ref{th:Hpers}, for all $z=(\rho, \theta, i) \in \mathcal{X}_+$
$$ \| \delta_z \tilde{P}_t - \Pi \circ \Psi^{-1} \|_{TV} \leq c ( 1 + \tilde{W}(z)) e^{- \kappa t}.$$
Now, for all $A \in \mathcal{B}( M \times E)$ and all $(x,i) \in M^* \times E$, $\delta_{x,i}P_t(A) - \Pi(A) = \delta_{\Psi(x,i)}\tilde{P}_t( \Psi(A)) - \Pi \circ \Psi^{-1}(\Psi(A))$, so that
\begin{align*}
\| \delta_{x,i} P_t - \Pi  \|_{TV} & = \| \delta_{\Psi(x,i)} \tilde{P}_t - \Pi \circ \Psi^{-1} \|_{TV}\\
& \leq c ( 1 + \tilde{W}(\Psi(x,i))) e^{- \kappa t}\\
& = c ( 1 + W(x,i)) e^{- \kappa t}.
\end{align*}
Then Theorem~\ref{th:persist3} is proved.
\qed

\prf \textbf{of Theorem~\ref{th:noncompactglobal}}

It suffices to show that Theorems \ref{th:Hpers} and \ref{th:Hnonpers} remain valid under  assumptions \ref{boundedjump} and \ref{Lyapunovinfinity}. For Theorem \ref{th:Hpers}, we show that Hypothesis 3 in \cite{Ben14} holds. That is, we have to check that there exists a continuous function $W : M \times E \to \RR_+$ with $\lim_{\|x\| \to \infty} W(x,i)=\infty$, a continuous function $LW : M \times E \to \RR_+$, $\alpha > 0$ and $C \geq 0$ such that 
\begin{enumerate}
\item[\textbf{(i)}] For every compact set $K \subset M$, there exists $W_K \in \D^2$ such that
\begin{enumerate}
\item[\textbf{(a)}] $W|_K = W_K|_K$ and ${\cal L}W_K|_K = LW|_K$,
\item[\textbf{(b)}] For all $x \in M$, $\sup \{ P_t( \Gamma W_K), \quad t \geq 0, \quad K \quad \text{compact} \} < \infty$
\end{enumerate}
\item[\textbf{(ii)}] $$ L W \leq - \alpha W + C. $$
\end{enumerate}
The only difference with Hypothesis \ref{Lyapunovinfinity} is that here $W_K$ has to be in $\D^2$. so we are done if we prove that $C_c^1 \subset \D^2$, which is equivalent to  $C_c^1 \subset \D$. Here, we use the weaker notion of domain given in \cite{Ben14} : a function $f$ is in $\D$ if :

\begin{enumerate}
\item $\mathcal{L} f(x,i) = \lim_{t \to 0} \frac{P_t f(x,i) - f(x,i)}{t}$ exists for all $(x,i) \in M \times E$;
\item  $\mathcal{L} f$ is continuous bounded;
\item $\sup_{0 < t \leq 1} \frac{1}{t} \| P_t f - f \| < \infty.$
\end{enumerate}
Let $f \in C_c^1$. Since the jumps rates are bounded, the proof of \cite[Proposition 2.1]{BMZIHP} adapts verbatim to the noncompact case provided the derivative of $f$ vanish outside a compact set - which is the case by definition of $C_c^1$. For Theorem \ref{th:Hnonpers}, we note that point (i) and (ii) in its proof are still valid since in our case, the set $\XX_0$ is compact (see Propositions 8.2 and 8.3 in \cite{Ben14}). Thus, point \textbf{(i)} of Theorem \ref{th:Hnonpers} can be shown by the same argument even if $\XX$ is not compact. Now, the existence of a Lyapunov function implies that there exists a compact set $K \subset M$ containing 0, such that, for all $(x,i) \in M \times E$, $\mathbb{P}_{(x,i)}( T_K <\infty) = 1$, where $T_K$ is the hitting time of $K$. Moreover, due to the accessibility of 0, for all neighbourhood $U$ of 0, there exists $\delta > 0$ such that, for all $(x,i) \in K \times E$, $\mathbb{P}_{(x,i)}( T_U <\infty) \geq \delta$. Hence, by Markov property, $\mathbb{P}_{(x,i)}( T_U <\infty) \geq \delta$ for all $(x,i) \in M \times E$ and point (ii) of Theorem \ref{th:Hnonpers} follows.
\qed
\section{Proof of Theorem \ref{th:cvexpo}}
\label{sec:proofcvexpo}
Before proving our convergence theorem, we first recall the definition of the Birkhoff part metric and some properties of monotone and subhomogeneous random dynamical systems given in the book of Chueshov~\cite{chueshov}. Let $D$ be a non-empty subset of $\{1,\ldots,d\}$ and let $\RR^d_{++,D}$ be the subset of $x \in \RR^d_+$ such that $x_i > 0$ if $i \in D$ and $x_i= 0$ otherwise. Then $\RR^d_{++,D}$ is called a \textit{part}. The \textit{Birkhoff part metric} is defined, for all $x, y  \in \RR^d_+$ by :
$$ p(x,y) = \max_{i \in D} |\log (x_i) - \log(y_i) |$$
if $x$ and $y$ are both in the same part $\RR^d_{++,D}$ for some $D$, and $p(x,y)= +\infty$ otherwise.
By monotony and strong subhomogeneity of $\Psi$,~\cite[Lemma 4.2.1]{chueshov} ensures that $\Psi$ is nonexpansive under the part metric on every part and strictly nonexpansive on $\RR^d_{++}$. In other words, for all $t \geq 0$, for all $\omega \in \Omega$, for all $D \subset \{1,\ldots,d\}$, for all $x, y \in \RR^d_{++,D}$,
$$ p(\Psi(t,\omega,x),\Psi(t,\omega,y)) \leq p(x,y),$$
and the inequality is strict if $D= \{1,\ldots,d\}$, $x \neq y$ and $t >0$. We would like to have a \textit{contraction}, meaning that there exist $\alpha \in (0,1)$ such that $ p(\Psi(t,\omega,x),\Psi(t,\omega,y)) \leq \alpha p(x,y)$. The following crucial lemma states that this is true if we restrain ourselves to compact subset of $\RR_{++}^d$.
\blem
Let $\varphi : \RR^d_{+} \to \RR^d_{+}$ be a $C^2$ monotone strongly subhomogeneous map and $K$ be a compact subset contained in $\RR_{++}^d$. Then $\varphi$ is a contraction for $p$ on $K$, that is :
$$ \tau_{K}(\varphi) := \sup_{x,y \in K, x\neq y} \frac{p(\varphi(x),\varphi(y))}{p(x,y)} < 1.$$
\label{lem:contractp}
\elem
\prf
First note that for all $x, y \in K$, with $x \neq y$, one has $\frac{p(\varphi(x),\varphi(y))}{p(x,y)} < 1$. In particular, by continuity of $p$ and $\varphi$, for all $\varepsilon >0$ there exists $\alpha < 1$ such that
\beq
\label{inlemsup1}
\sup_{x,y \in \Delta_{\varepsilon}(K)} \frac{p(\varphi(x),\varphi(y))}{p(x,y)} \leq \alpha, \eeq
where $\Delta_{\varepsilon}(K) = \{ (x,y) \in K^2 \: : \: p(x,y) \geq \varepsilon \}$ is compact. It remains to prove that such a bound holds when $x$ and $y$ are close, uniformly in $x \in K$.
To do so, we use the following fact: a monotone map $\varphi$ is strongly sublinear if and only if, for all $x \gg 0$, $D \varphi(x) x \ll \varphi(x)$ (see e.g \cite[Proposition 4.1.1]{chueshov} or \cite[Proposition 6]{BS09}). Componentwise, this means that for all $i$,
\beq
\label{inlemcomp}
\frac{\langle \nabla \varphi_i(x), x \rangle}{\varphi_i(x)} < 1.
\eeq
By Taylor expansion, for all $i$ and all $x, y \in K$,
$$ \log \varphi_i(y) - \log \varphi_i(x) = \frac{\langle \nabla \varphi_i(x), y- x \rangle}{\varphi_i(x)} + R_i(x,y) \| x -y \|^2,$$
where $R_i$ is continuous, thus uniformly bounded on $K^2$ by some constant $C$.

Moreover, one can easily check that for all $\frac{1}{2M} \leq u \leq 2M$, one has
$$ | u -1 | \leq e^{|\log u|} - 1 \leq | \log u | ( 1 + M  | \log u |).$$ 
Now there exists $M$ such that for all $x, y \in K$ and $k$, $\frac{1}{2M} \leq y_k/x_k \leq 2M$. Thus, for all $k$,
\beq
\label{inlem}
|y_k - x_k | \leq x_k (1+ M p(x,y))p(x,y).
\eeq
For all $x,y \in \RR^d_{++}$ and $x \neq y$, there exists $i$ such that
\begin{align*}
 \frac{p(\varphi(x),\varphi(y))}{p(x,y)} & = \frac{\vert \frac{\langle \nabla \varphi_i(x), y- x \rangle}{\varphi_i(x)} + R_i(x,y) \| x -y \|^2 \vert}{p(x,y)}\\
 & \leq \frac{\vert \langle \nabla \varphi_i(x), y- x \rangle \vert}{\varphi_i(x) p(x,y)} + \vert  R_i(x,y) \vert \frac{\| x -y \|^2 }{p(x,y)}
 \end{align*}
Now by \eqref{inlem} and nonnegativity of $\nabla \varphi_i(x)$ (recall $\varphi$ is monotone), we have for all $x, y \in K$, for all  $x \neq y$,
$$ \frac{p(\varphi(x),\varphi(y))}{p(x,y)} \leq \frac{\langle \nabla \varphi_i(x), x (1+M p(x,y)) \rangle}{\varphi_i(x)} + C \frac{\| x -y \|^2 }{p(x,y)}.$$
Inequality \eqref{inlemcomp}, continuity of $\varphi$ and compactness of $K$ imply that there exists a constant $\tau <1$ such that, for all $x \in K$ and all $i$,
$$\frac{\langle \nabla \varphi_i(x), x \rangle}{\varphi_i(x)} \leq \tau, $$
and thus
$$ \frac{p(\varphi(x),\varphi(y))}{p(x,y)} \leq \tau (1+ M p(x,y)) + C \frac{\| x -y \|^2 }{p(x,y)}.$$
By compactness of $K$, $p(x,y)$ and $\frac{\| x -y \|^2 }{p(x,y)}$ converges to 0 uniformly in $x \in K$ when $y$ converges to $x$. Thus, we can find $\varepsilon > 0$ such that $ \tau' = \sup_{x \in K, y \in  B_{K}(x,\varepsilon) \setminus \{ x \}}  \tau (1+M p(x,y)) + C \frac{\| x -y \|^2 }{p(x,y)} < 1$, where $B_{K}(x,\varepsilon)$ is the intersection of the ball of center $x$ and radius $\epsilon$ with $K$. In other words,
\beq
\label{inlemsup2}
\sup_{x,y \in \Delta_{\varepsilon}^{\mathrm{c}}(K)} \frac{p(\varphi(x),\varphi(y))}{p(x,y)} \leq \tau'.
\eeq
Combining \eqref{inlemsup1} and \eqref{inlemsup2} gives the result with $\tau_K(\varphi) = \max ( \alpha, \tau') < 1$.
\qed

Recall that $\mathcal{Y}= [0,1]^d \setminus \{0\} \times E$ and set $d : \mathcal{Y}^2 \to [0,1]$ the distance defined by
$$d((x,i),(y,j)) = \1_{i \neq j} + \1_{i=j}(\frac{p(x,y)}{C} \wedge 1),$$
where $C$ is a constant to be chosen later  and $p(x,y)$ is the Birkhoff part metric. Define also $V : \mathcal{Y} \to \RR_+$ with $V(x,i) = \|x\|^{- \theta}$ where $\theta$ is given in Theorem~\ref{th:persist1} and the function $\tilde{d} :  \mathcal{Y}^2 \to \RR_+$ by
 $$\tilde{d}(z,\tilde{z}) = \sqrt{d(z,\tilde{z})(1+V(z)+V(\tilde{z}))}.$$
As already mentioned, Theorem \ref{th:cvexpo} is a consequence of the weak form of Harris' theorem due to Hairer, Mattingly and Scheutzow~\cite[Theorem 4.8 and remark 4.10]{hms}. More precisely, it states that point \textbf{(i)} of  Theorem \ref{th:cvexpo} holds, provided the three following assumptions are verified (here we let $P_t$ denoted $P_t^Z$) :
\bdes
\item[A1] V is a \textit{Lyapunov function} for $P_t$, that is there exists $C_V ,\gamma , K_V, t_0 >0$ such that for all $t \geq t_0$, for all $z \in \mathcal{X}$,
$$ P_t V(z) \leq C_V e^{- \gamma t} V(x) + K_V;$$
\item[A2] There exists $t^*>t_*>0$ such that for all $t \in [t_*,t^*]$, the level set $ A_V = \{ z \in  \mathcal{X} \: : \: V(x) \leq 4 K_V \}$ are \textit{$d$-small } for $P_t$, meaning that there exists $\varepsilon >0$ such that for all $z, \tilde{z} \in A_V$,
$$ \mathcal{W}_d(\delta_z P_t, \delta_{\tilde{z}}P_t) \leq 1 - \varepsilon;$$
\item[A3]  For all $t \in [t_*,t^*]$, $P_t$ is \textit{contracting} on $A_V$, meaning that there exists $\alpha \in (0,1)$ such that for all $z, \tilde{z} \in A_V$ with $d(z, \tilde{z}) < 1$,
$$  \mathcal{W}_d(\delta_z P_t, \delta_{\tilde{z}}P_t) \leq \alpha d(z, \tilde{z}).$$
Moreover, $P_t$ is nonexpansive on $\mathcal{X}$, that is for all $z, \tilde{z} \in \mathcal{X}$,
$$  \mathcal{W}_d(\delta_z P_t, \delta_{\tilde{z}}P_t) \leq  d(z, \tilde{z}).$$
 \edes
\brem {\rm
In~\cite[Theorem 4.8]{hms}, the hypothesis \textbf{A1} and \textbf{A3} are a little bit stronger : \textbf{A1} should holds for every $t \geq 0$, and the contraction in \textbf{A3} should holds on the whole space $\mathcal{X}$ for $d(z, \tilde{z}) < 1$. However, a quick look at the proof given in~\cite{hms} shows that it is enough to have the Lyapunov function for $t$ large, and that when $z,  \tilde{z} $ are such that $1+V(z)+V(\tilde{z}) \geq 4K_V$, the proof "Far from the origin" is true independently from the fact that $d(z, \tilde{z}) < 1$ or $d(z, \tilde{z}) \geq 1$
}
\erem

To prove Theorem~\ref{th:cvexpo} it is thus sufficient to show that  \textbf{A1} to \textbf{A3} are satisfied. For  \textbf{A1}, it is a consequence of a stochastic persistence lemma. For  \textbf{A2}, we show that a good choice of the constant $C$ appearing in the definition of $d$ is sufficient to have the small set. Finally,  \textbf{A3} is a consequence of the contracting properties of $\Psi(t,\omega)$.

\bigskip
\noindent
\prf  \textbf{of Theorem~\ref{th:cvexpo}}
\paragraph*{A1}
We have the following lemma :
\blem
For $0< \alpha < \lambda_1$,  there exists $T > 0$, $\eps >0$ and $C > 0$ such that, for all $t  \in [T, 3T/2]$, for all $z \in \mathcal{Y}^{\varepsilon}_0$,
$$ P_t V(z) \leq e^{ \theta t ( \frac{t}{T}-1) \alpha} V(z),$$
where $\theta = \frac{\alpha}{C T}$, $\mathcal{Y}^{\varepsilon}_0= \{ (x,i) \in \mathcal{Y} \: : \: \|x\| < \varepsilon \}$ and $V(x,i) = \|x\|^{- \theta}$.
\elem

\prf
Follows the lines of the proof given in \cite[Lemma 3.5]{BL16}.
\qed

In particular, putting $\gamma = \frac{\theta \alpha}{4}$, then for all $t \in [T, 3T/2]$, for all $z \in  \mathcal{Y}^{\varepsilon}_0$,
$$ P_t V(z) \leq e^{ \gamma t } V(z).$$
Now by Feller continuity of $P_t$ and compactness of $[T, 3T/2] \times \mathcal{Y} \setminus \mathcal{Y}^{\varepsilon}_0$
$$ \tilde{C} = \sup_{(t,z) \in [T, 3T/2] \times \mathcal{Y} \setminus \mathcal{Y}_{\varepsilon}} P_t V(z) - V(z) < \infty,$$
and, for all $t \in [T, 3T/2]$ and all $z \in \mathcal{Y}$,
$$ P_t V(z) \leq e^{ \gamma t } V(z) + \tilde{C}.$$
If $t \geq 2T$, then there exists $s \in [T, 3T/2]$ and $n \geq 1$ such that $t = ns$. Thus
$$ P_t V(z) = P_{ns} V(z) \leq e^{ \gamma n s } V(z) + \sum_{k=0}^{n-1} e^{ \gamma k s } \tilde{C},$$
proving \textbf{A1} with $t_0 = 2T$ and $K_V = \frac{1}{1-e^{-\gamma T}} \tilde{C}$.
\paragraph*{A2}
Set $M_V = \{ x \in [0,1]^d \setminus \{0\} \: : \: \|x\|^{-\theta} \leq 4 K_V \}$. We first prove that for all $t^*>t_* > 0$, there exists a compact set contained in $\RR_{++}^d$ such that for all $t \in [t_*,t^*]$, and all $\omega \in \Omega$, $\Psi(t,\omega, M_V)$ is included in this compact. For this, let $S_{M_V}$ denotes the set of all the solutions of the differential inclusion
$$\left \{\begin{array}{l}
\dot{\eta}(t) \subset \mathsf{co}(\tilde{\mathrm{F}})(\eta(t)) \\
\eta(0) = x,\end{array} \right.$$
with $x \in M_V$. Then because $M_V$ is compact, $S_{M_V}$ is a non avoid compact subset of $\mathcal{C}(\RR_+,\RR^d)$ (see e.g Aubin and Cellina \cite[Section 2.2 Theorem 1]{aubincellina}). This implies that $\mathbf{\Psi}_{[t_*,t^*]}(M_V) = \{ \eta_t \: : \: t \in [t_*,t^*],  \eta \in S_{M_V} \}$ is a compact set of $[0,1]^d$. Moreover, by strong monotony of $\eta_t$,  $\mathbf{\Psi}_{[t_*,t^*]}(M_V)$ is included in $(0,1]^d$ and for all $t \in [t_*,t^*]$, $\omega \in \Omega$,  $\Psi(t,\omega, M_v) \subset \mathbf{\Psi}_{[t_*,t^*]}(M_V)$.
Now by compactness of  $\mathbf{\Psi}_{[t_*,t^*]}(M_V)$ and continuity of $p$, there exist $K >0$ such that for all $t \in [t_*,t^*]$,
\beq
\label{eq:inegK}
\sup_{x,y \in M_v; \omega, \omega' \in \Omega} p( \Psi(t,\omega,x), \Psi(t,\omega',y)) \leq \sup_{a,b \in \mathbf{\Psi}_{[t_*,t^*]}(M_V)} p(a,b) = K.
\eeq
To prove \textbf{A2}, for any $(z,\tilde{z})=((x,i),(y,j)) \in \mathcal{Y}^2$, we consider the coupling $(Z_t,\tilde{Z}_t) = ((X_t,I_t),(Y_t,J_t))$ of $\delta_z P_t$ and $\delta_{\tilde{z}} P_t$ construct as follows. If $i = j$, then $I_t = J_t$ for all $t\geq 0$. If $i \neq j$, then $I_t$ and $J_t$ evolves independently until the first meeting time $T$ and then are stick together for ever. In other words,
$$ \mathbb{P}_{i,j} ( I_t \neq J_t) =  \mathbb{P}_{i,j}( T > t).$$
This is the coupling considered in \cite{ecp}. As stated in~\cite[Lemma 2.1]{ecp}, we easily control the above probability : there exists $\rho > 0$ such that for all $i,j \in E$ and all $t \geq 0$,
$$ \mathbb{P}_{i,j} ( I_t \neq J_t) =  \mathbb{P}_{i,j}( T > t) \leq e^{-\rho t}.$$
Let $(z,\tilde{z})=((x,i),(y,j)) \in A_V^2$ and $t \in [t_*,t^*]$. Then
\begin{align*}
\mathcal{W}_d(\delta_z P_t, \delta_{\tilde{z}} P_t) & \leq \mathbb{E}_{(z,\tilde{z})}(d( Z_t, \tilde{Z}_t))\\
& \leq \mathbb{P}_{i,j} ( I_t \neq J_t) + \mathbb{E}_{(z,\tilde{z})}(\frac{p(X_t,Y_t)}{C})\\
& \leq e^{- \rho t} + \frac{K}{C},
\end{align*}
where the last inequality comes from~\eqref{eq:inegK}. Thus, choosing $C = \frac{K}{1-2e^{-\rho t_*}}$, one has
$$ \mathcal{W}_d(\delta_z P_t, \delta_{\tilde{z}} P_t) \leq 1 + e^{-\rho t} - 2 e^{-\rho t_*} \leq 1- e^{-\rho t_*}, $$
proving \textbf{A2} with $\varepsilon =  e^{-\rho t_*}$.
\paragraph*{A3} We first prove that $P_t$ is nonexpansive on $\mathcal{Y}$. Is suffices to show the result for $(z,\tilde{z})$ such that $d(z,\tilde{z}) < 1$, the bound being trivial otherwise. In particular, $i=j$ where $z=(x,i)$ and $\tilde{z}=(y,j)$, and $d(z,\tilde{z})= \frac{p(x,y)}{C}<1$, which implies that $x$ and $y$ are in the same part. We consider the same coupling $(Z_t,\tilde{Z}_t)$ as above.  Then because $i=j$, $I_t = J_t$ and thus $X_t = \Psi(t,\omega,x)$ and $Y_t = \Psi(t,\omega,y)$, and so by nonexpansivity of $\Psi(t,\omega)$ on every part, one has $p(\Psi(t,\omega,x), \Psi(t,\omega,y)) \leq p(x,y)$, which gives the result for $P_t$.

Now we prove that $P_t$ is a contraction on $A_V$. Let  $t \in [t_*,t^*]$ and $(z,\tilde{z})\in A_V^2$ such that $d(z,\tilde{z}) < 1$. In addition with the consequences cited above, this also implies that $x,y \in M_V$. Choose $0<t_0 < t_*$, then one has
\begin{align*}
p(\Psi(t,\omega,x), \Psi(t,\omega,y)) & = p (\Psi(t-t_0+t_0,\omega,x), \Psi(t-t_0+t_0,\omega,y))\\
& \leq p(\Psi(t-t_0, \Theta_{t_0}\omega)\Psi(t_0,\omega,x), \Psi(t-t_0, \Theta_{t_0}\omega)\Psi(t_0,\omega,y))\\
& \leq \tau_{\mathbf{\Psi}_{t_0}(M_V)}(\Psi(t-t_0, \Theta_{t_0}\omega)) p(\Psi(t_0,\omega,x),\Psi(t_0,\omega,y))\\
& \leq \tau_{\mathbf{\Psi}_{t_0}(M_V)}(\Psi(t-t_0, \Theta_{t_0}\omega)) p(x,y),
\end{align*}
where $\tau_{\mathbf{\Psi}_{t_0}(M_V)} (\Psi(t-t_0, \Theta_{t_0}\omega)) < 1$ is the contraction constant given by Lemma~\ref{lem:contractp} on the compact  $\mathbf{\Psi}_{t_0}(M_V) \subset \RR^d_{++}$. Because $\tau_{\mathbf{\Psi}_{t_0}(M_V)} (\Psi(t-t_0, \Theta_{t_0}\omega)) < 1$ for every $\omega$, then
$$\alpha = \max_{i} \mathbb{E}_{i}[\tau_{\mathbf{\Psi}_{t_0}(M_V)} (\Psi(t-t_0, \Theta_{t_0}\omega))] < 1,$$
and
$$ \mathcal{W}_d(\delta_z P_t,\delta_{\tilde{z}} P_t) \leq \mathbb{E}_{(x,i),(y,j)} ( \frac{p(\Psi(t,\omega,x), \Psi(t,\omega,y))}{C} \leq \alpha \frac{p(x,y)}{C}= \alpha d(z,\tilde{z}), $$
proving \textbf{A3} and the \textbf{(i)} of the theorem.

Because $\lambda_1 >0$, Theorem \ref{th:persist1} insures existence of an invariant measure for $P_t$ on $\mathcal{Y}$. The uniqueness of the invariant measure and thus point \textbf{(ii)} follows immediately from point \textbf{(i)}.
\qed
\section{Appendix}
\mylabel{sec:appendix}
\subsection{Proof of Proposition~\ref{th:Perron}}
 Recall (see section \ref{sec:epidemic}) that $\RR^d_{++}$ denotes the interior of $\RR^d_+$,
(i.e~ the cone of positive vectors). Set $S^{d-1}_+ = S^{d-1} \cap \RR^d_+$ and $S^{d-1}_{++} = S^{d-1} \cap \RR^d_{++}$. The principal tool is the {\em projective} or \textit{Hilbert metric}  $d_H$ on $ \RR^d_{++}$ (see Seneta \cite{Seneta}) defined by
\[ d_H(x,y) = \log \frac{\max_{1 \leq i \leq d} x_i/y_i}{\min_{1 \leq i \leq d} x_i/y_i}. \]
Note that
\begin{equation}
d_H(\frac{x}{\|x\|},\frac{y}{\|y\|}) = d_H(x,y)
\label{hilbproj}
\end{equation}
so that $d_H$ is not a distance on $\RR^d_{++}.$ However its restriction to $S^{d-1}_{++}$ is. Furthermore, for all $x$, $y \in S^{d-1}_{++}$,
\begin{equation}
 \| x - y \| \leq \mathrm{e}^{d_H(x,y)} - 1.
 \label{hilbnorm}
\end{equation}
Let ${\cal M}_{+}$ denote the set of $d \times d$ Metzler matrices having   {\bf positive} diagonal entries,
 and  let ${\cal M}_{++} \subset  {\cal M}_{+}$ denote  the set of matrices having positive entries.
 By a theorem of Garret Birkhoff,  there exists a continuous map $\tau :  {\cal M}_{++} \mapsto ]0,1[$ such that for all $T \in {\cal M}_{++},$ and all  $x, y \in \RR^d_{++}$
  \beq
  \label{contractdH}
   d_H(Tx, Ty) \leq \tau[T] d_H(x,y)
   \eeq
  The number $\tau[T]$ is usually called  the \textit{Birkhoff's contraction coefficient} of $T,$ and is given by an explicit formulae
   (see e.g \cite{Seneta}, Section 3.4) which is unneeded here.

We extend $\tau$ to a measurable map $\tau : {\cal M}_{+} \mapsto ]0,1]$  by setting $\tau[T] = 1$ for all  $T \in {\cal M}_{+} \setminus {\cal M}_{++}.$  By density of  ${\cal M}_{++}$ in ${\cal M}_{+}$ and continuity of $d_H$ on $\RR^d_{++}$  it is easy to see that (\ref{contractdH}) extends to ${\cal M}_{+}.$

For each $\omega \in \Omega,$ the map  $t \mapsto \varphi(t,\omega)$ is solution to the matrix valued differential equation
\beq
\label{matrixpdmp}
\forall t \geq 0, \, \frac{dM}{dt} = A^{\omega_t} M, M_0 = I_d.
\eeq
Thus, $$\varphi(t,\omega) \in {\cal M}_{+}$$ for all $t \geq 0.$ Indeed, for all $i \in E$ and $r > 0$ large enough  $A^i+ r I_d \in {\cal M}_{+},$ so that $e^{tA^i} =  e^{-rt} e^{t (A^i + r I_d)} \in {\cal M}_{+}.$

We claim that there exists a Borel set $\tilde{\Omega} \subset \Omega$ with  $\P^J_i(\tilde{\Omega}) = 1$ for all $i \in E,$ and  such that  for all $\omega \in \tilde{\Omega} :$
\bdes
\iti $\exists n \in \NN \,  \varphi(n, \omega) \in {\cal M}_{++};$
\itii $\forall n \in \NN \,  \limsup_{t \rar \infty} \displaystyle{ \frac{\log{\tau[\varphi(t,\mathbf{\Theta}_{n}(\omega))]}}{t}} < 0.$
\edes
Before proving these assertions  let us show how they imply the result to be proved.
For all $\omega \in \tilde{\Omega}$  and $n$   given by $(i),$
$$\varphi(t+n,\omega)  = \varphi(t, \mathbf{\Theta}_{n}(\omega)) \varphi(n,\omega)  \in {\cal M}_{++}$$
as the product of an element of ${\cal M}_{+}$ with an element  of ${\cal M}_{++}.$
 Thus, by $(ii)$, for all $\omega \in \tilde{\Omega}$ and $x, y \in \RR^d_+ \setminus \{0\}$
\begin{equation}
\label{cvtheta}
\limsup_{t \rar \infty} \frac{1}{t} \log  d_H( \varphi(t + n,\omega)x , \varphi(t + n,\omega)y ) < 0.
\end{equation}
For $x \in S^{d-1}_+$ set $$\Phi(t,\omega)x = \frac{\varphi(t,\omega) x}{\|\varphi(t,\omega) x\|}.$$
 Let $f : S^{d-1}_+ \times E \to \RR$ be a continuous map. It follows from \eqref{cvtheta}, ~\eqref{hilbproj}, ~\eqref{hilbnorm} and the continuity of $f$ that  $$| f(\Phi(t,\omega)x, \omega_t) -  f(\Phi(t,\omega)y, \omega_t)| \to 0$$  for all $x,y \in S^{d-1}_+$ and $\omega \in \tilde{\Omega}.$ Moreover, $$P_t^{(\Theta,J)}f(x,i) = \mathbb{E}_{i}^J( f(\Phi(t,\omega)x, \omega_t)),$$ and thus $$\lim_{t \rar \infty}  P_t^{(\Theta,J)}f(x,i) - P_t^{(\Theta,J)}f(y,i) = \lim_{t \rar \infty} \mathbb{E}_{i}^J( f(\Phi(t,\omega)x, \omega_t) - f(\Phi(t,\omega)y, \omega_t)) =  0$$ by dominated convergence. Now take $\mu$, $\nu \in {\cal P}_{inv}^{(\Theta, J)}$. Then one has
\begin{equation}
\label{cvpt}
 \lim_{t \to \infty} \sum_i p_i  \int_{ (S^{d-1}_+)^2}  \left( P_t^{(\Theta,J)}f(x,i) - P_t^{(\Theta,J)}f(y,i) \right) \mu(\mathrm{d}x | i) \nu(\mathrm{d}y | i) = 0,
 \end{equation}
where $\mu(\cdot| i) = \mu^i( \cdot)/p_i$. But by invariance of $\mu$ and $\nu$, the left-hand side of~\eqref{cvpt} equals $\mu f - \nu f$ for all $t$, giving $\mu f = \nu f$ for all continuous $f.$  This proves unique ergodicity of $(\Theta,J).$

We now pass to the proofs of assertions $(i)$ and $(ii)$ claimed above.

Irreducibility of $\overline{A}$ implies that $e^{\overline{A}} \in {\cal M}_{++}.$ Let  $\mathcal{U} \subset {\cal M}_{++}$  be a compact neighborhood of $e^{\overline{A}}.$
Since $\overline{A}.M \in \mathsf{co}(A^i)(M)$, it follows from
the Support Theorem (~\cite[Theorem 3.4]{BMZIHP}), applied to the PDMP (\ref{matrixpdmp}), that  for all $i \in E$
\[\P^J_i \{\omega \in \Omega : \: \varphi(1,\omega) \in \mathcal{U}\} > 0. \]
Thus, by  the Markov property or the conditional version of the Borel Cantelli Lemma, for  $\P^J_i$ almost all $\omega,$
$\varphi(1,\mathbf{\Theta}_n(\omega)) \in \mathcal{U}$ for infinitely many $n,$ and consequently, for $n$ large enough
$$\varphi(n,\omega) = \varphi(1, \mathbf{\Theta}_{n-1}  \omega) \ldots  \varphi(1, \omega) \in {\cal M}_{++}.$$ This proves assertion $(i).$ By the cocycle property and Birkhoff ergodic theorem, for $\P_p^J$ (hence $\P_i^J$) almost all $\omega$
$$\limsup_{t \rar \infty} \frac{1}{t} \log(\tau[\varphi(t,\omega)])  \leq   \limsup_{n \rar \infty} \frac{1}{n} \log(\tau[\varphi(n,\omega)]) \leq
 \limsup_{n \rar \infty} \frac{1}{n} \sum_{k = 1}^n \log \left (\tau[\varphi(1, \mathbf{\Theta}_{k-1}( \omega))]\right )$$
  $$ = \E_p^J(\log(\tau[\varphi(1,\omega)])) \leq \sup_{ M \in {\cal U}} \log (\tau[M])  \P_p^J(\omega  \in \Omega \: : \varphi(1,\omega) \in {\cal U}) < 0.$$ Replacing $\omega$ par $\mathbf{\Theta_n(\omega)}$ proves assertion $(ii).$
 \qed

\subsection{Proof of Lemma \ref{lem:moyenne}}
Before proving Lemma \ref{lem:moyenne}, we prove the following lemma, which is a consequence of results from Freidlin and Wentzell \cite{FW12}.

\blem
Assume the switching rates are constant and depend on a small parameter $\varepsilon$ : $a_{i,j}^{\varepsilon} = a_{i,j}/{\varepsilon} $ where $(a_{i,j})$ is an irreducible matrix with invariant probability $p$. Denote by $(X^{\varepsilon},J^{\varepsilon})$ the PDMP associated with $a_{i,j}^{\varepsilon}$ given by \eqref{eq:pdmp}. Let $\Psi$ denote the flow induced by the average vector field $F^p := \sum_i p_i F^i$
Then for all $\delta > 0$ and all $T > 0$,

\beq
 \lim_{ \varepsilon \to 0} \mathbb{P}_{(x,i)} \left( \max_{ 0 \leq t \leq T} | X_t^{\varepsilon} - \Psi_t(x) | > \delta \right) = 0,
 \eeq
uniformly in $(x,i) \in M \times E$.
\label{lem:average}
\elem

\prf
According to \cite[Chapter 2 Theorem 1.3]{FW12}, it suffices to show that for all $\delta > 0$ and all $T > 0$,

\beq
 \lim_{ \varepsilon \to 0} \mathbb{P}^{J}_{i} \left( \left| \int_{t_0}^{t_0+T}(F^{J_t^{\varepsilon}}(x) - F^p(x)) \mathrm{d}t \right| > \delta \right) = 0,
 \label{eq:cvFW}
 \eeq
uniformly in $t_0 >0$ and $(x,i) \in M \times E$. Note that

\begin{align*}
\left| \int_{t_0}^{t_0+T}(F^{J_t^{\varepsilon}}(x) - F^p(x)) \mathrm{d}t \right| & = \left| \int_{t_0}^{t_0+T}(\sum_j F^{j}(x) \1_{J_t^{\varepsilon}=j} - \sum_ j p_j F^j(x)) \mathrm{d}t \right|\\
& \leq \sum_j \| F^j \|_{\infty}  \left| \int_{t_0}^{t_0+T} (\1_{J_t^{\varepsilon}=j} - p_j) \mathrm{d}t \right|,
\end{align*}
so \eqref{eq:cvFW} is proven if we show that $\int_{t_0}^{t_0+T} \1_{J_t^{\varepsilon}=j} \mathrm{d}t$ converges in probability to $p_j T$ uniformly in $t_0 > 0$. By Fubini's Theorem and invariance of $p$, $ \mathbb{E}^{J}_p \left( \int_{t_0}^{t_0+T} \1_{J_t^{\varepsilon}=j} \mathrm{d}t \right) = p_j T$, so Bienaymé - Tschebischev inequality gives

$$  \mathbb{P}^{J}_{i} \left( \left| \int_{t_0}^{t_0+T} (\1_{J_t^{\varepsilon}=j} - p_j) \mathrm{d}t \right| > \delta \right) \leq \frac{V^J_p(\int_{t_0}^{t_0+T} (\1_{J_t^{\varepsilon}=j}\mathrm{d}t)}{\delta},$$
where $V_p^J$ is the variance associated to $\mathbb{E}_p^J$.  Hence we can conclude if $\mathbb{E}^{J}_p \left[ \left( \int_{t_0}^{t_0+T} \1_{J_t^{\varepsilon}=j} \mathrm{d}t \right)^2 \right]$ converges to $(p_j T)^2$ uniformly in $t_0 > 0$.

Denote by $Q$ the intensity matrix of $J^{1}$, then for all $\varepsilon >0$, the intensity matrix of $J^{\varepsilon}$ is $Q/\varepsilon$ and for all $i,j \in E$ and $t \geq 0$,

$$ \mathbb{P}_{i}(J^{\varepsilon}_t = j ) = \left( e^{\frac{t}{\varepsilon}Q} \right)_{i,j}.$$
By ergodicity of $ J^{\varepsilon}_t$, the above quantity goes to $p_j$ when $t \to \infty$ so also for every fixed $t$ when $\varepsilon$ goes to 0. Now we have

\begin{align*}
{E}^{J}_p \left[ \left( \int_{t_0}^{t_0+T} \1_{J_t^{\varepsilon}=j} \mathrm{d}t \right)^2 \right] & = 2 \int_{t_0}^{t_0+T} \int_{t_0}^{t} \mathbb{P}_p \left( J^{\varepsilon}_u = j; J^{\varepsilon}_t = j \right) \mathrm{d}u \mathrm{d}t\\
& =  2 \int_{t_0}^{t_0+T} \int_{t_0}^{t} \mathbb{P}_j \left( J^{\varepsilon}_{t-u} = j \right) p_j \mathrm{d}u \mathrm{d}t\\
& =  2 \int_{t_0}^{t_0+T} \int_{t_0}^{t} \left( e^{\frac{t-u}{\varepsilon}Q} \right)_{j,j}  p_j \mathrm{d}u \mathrm{d}t,
\end{align*}
where the second inequality resulted from the Markov property. Now because for all $t_0$, $t-u \in [0,T]$, $\left( e^{\frac{t-u}{\varepsilon}Q} \right)_{j,j} $ converges almost everywhere to $p_j$  and thus the lemma is proven by dominated convergence.
\qed

With the notation of the preceding lemma,
let $$\mu^{\varepsilon} \in {\cal P}_{inv}^{(X^{\varepsilon},J^{\varepsilon})}, \, \nu^{\varepsilon} = \sum_i \mu^{i,\varepsilon}.$$
The proof of the next lemma is similar to the proof of \cite[Corollary 3.2]{B98recu}.
\blem
Let  $\nu$ a limit point of $(\nu^{\varepsilon})$ when $\varepsilon \rar 0.$ Then
 $\nu$  is an invariant measure of $F^p$.
\elem
\prf
 For notational convenience, we assume that $\nu^{\varepsilon}$ converges to $\nu$.

Let $g : M \to \RR$ be a continuous map, then for all $t > 0$ and all $\varepsilon >0$,
\begin{align*}
\left| \int g(\Psi_t) d \nu - \int g d \nu \right| &
\leq  \left| \int g(\Psi_t) d \nu - \int g  d \nu^{\varepsilon} \right| + \left| \int g d \nu^{\varepsilon} - \int g d \nu \right|\\
& \leq \left| \int g(\Psi_t) d \nu - \int g(\Psi_t)  d \nu^{\varepsilon} \right| +
\left| \int g(\Psi_t)  d \nu^{\varepsilon} - \int \mathbb{E}(g (\Theta_t^{\varepsilon}) ) d \nu^{\varepsilon} \right|\\
& \quad + \left| \int g d \nu^{\varepsilon} - \int g d \nu \right|,
\end{align*}
where we have use invariance of $\nu$ and $\nu^{\varepsilon} $. The first and the last term of the right hand side converge to 0 by definition of $\nu$, and the second one also converges to 0 by Lemma \ref{lem:average}.
\qed

Now let $\mu$ be a limit point of $(\mu^{\varepsilon})$. For notational convenience, we assume that $\mu^{\varepsilon}$ converges to $\mu$. We prove that  $\mu = \nu \otimes p $, which implies Lemma \ref{lem:moyenne}. For every continuous  $f : M \times E \to \RR$, every $t \geq 0$ and $\varepsilon >0$, one has
\begin{align*}
\mu^{\varepsilon} f - \mu f & = \int_{M \times E} \mathbb{E}_{(x,i)} \left( f_{J^{\varepsilon}_t}(X_t^{\varepsilon}) \right) d \mu^{\varepsilon}(x,i) - \sum_j p_j \int_{M} f_j(\Psi_t(x)) d \nu(x)\\
& = \int_{M \times E} \mathbb{E}_{(x,i)} \left( f_{J^{\varepsilon}_t}(X_t^{\varepsilon}) \right) d \mu^{\varepsilon}(x,i) - \int_{M \times E} \mathbb{E}_{(x,i)} \left( f_{J^{\varepsilon}_t}(\Psi_t) \right) d \mu^{\varepsilon}(x,i)\\
& \quad + \int_{M \times E} \mathbb{E}_{(x,i)} \left( f_{J^{\varepsilon}_t}(\Psi_t) \right) d \mu^{\varepsilon}(x,i) - \sum_j p_j \int_{M \times E} f_j(\Psi_t(x)) d \mu^{\varepsilon}(x,i)\\
& \quad + \sum_j p_j \int_{M \times E} f_j(\Psi_t(x)) d \mu^{\varepsilon}(x,i) - \sum_j p_j \int_{M} f_j(\Psi_t(x)) d \nu(x)  \\
& = A + B + C.
\end{align*}
We have
$$ \sup_{(x,i) \in M \times E} \mathbb{E}_{(x,i)} \left| f_{J^{\varepsilon}_t}(X_t^{\varepsilon}) -  f_{J^{\varepsilon}_t}(\Psi_t) \right)| \leq \max_j \sup_{(x,i) \in M \times E} \mathbb{E}_{(x,i)} \left( f_j(X_t^{\varepsilon}) -  f_j(\Psi_t) \right),$$
where the right hand side converges to 0 when $\varepsilon$ goes to 0 thanks to Lemma \ref{lem:average}, so $A$ converges to $0$. Next,
$$ |B | \leq  \sum_j \int_{M \times E} \left| \mathbb{P}_i ( J^{\varepsilon}_t = j) - p_j \right| | f_j( \Psi_t(x)) | d \mu^{\varepsilon}(x,i),$$
because $ \mathbb{E}_{(x,i)} \left( f_{J^{\varepsilon}_t}(\Psi_t) \right) = \sum_j  \mathbb{P}_i ( J^{\varepsilon}_t = j) f_j( \Psi_t(x))$. Thus $B$ converges to $0$ because $\left| \mathbb{P}_i ( J^{\varepsilon}_t = j) - p_j \right| $ converges to 0 uniformly in $i$ and $j$.
Finally, by definition of $\nu^{\varepsilon}$
$$ C =  \int_{M} \sum_j p_j f_j(\Psi_t(x)) d \mu^{1, \varepsilon}(x,i) -  \int_{M} \sum_j p_j f_j(\Psi_t(x)) d \nu(x),$$
proving that $C$ converges to 0 by definition of $\nu$ and thus the Lemma.
\qed

\section*{Acknowledgments}
This work was supported by the SNF grant $2000020 - 149871/1$. We thank Janusz Mierczy{\'n}ski for many valuable comments on section 2, especially proof of proposition \ref{th:lyapou}. We thank two anonymous referees for their useful comments. ES thanks Carl-Erik Gauthier for many discussions on this subject.

\bibliographystyle{amsplain}
\bibliographystyle{imsart-nameyear}
\bibliographystyle{nonumber}
\bibliography{RandomSwitch}

\end{document}